\documentclass[final]{amsart}
\usepackage{url,enumerate}

\newcommand{\eqinlaw}{\stackrel{\mathcal{L}}{=}}
\newcommand{\toinlaw}{\xrightarrow{\mathcal{L}}}
\newcommand{\E}[1]{\mathbf{E}\,{#1}}
\newcommand{\Var}[1]{\mathbf{Var}\,{#1}}
\newcommand{\mun}[1]{\mu_n^{[{#1}]}}
\newcommand{\tmun}[1]{\tilde{\mu}_n^{[{#1}]}}
\newcommand{\muz}[1]{\mu^{[{#1}]}(z)}

\newcommand{\tmuz}[1]{\tilde{\mu}^{[{#1}]}(z)}
\newcommand{\trz}[1]{\tilde{r}^{[{#1}]}(z)}

\newcommand{\A}{\mathcal{A}}
\newcommand{\N}{\mathcal{N}}
\newcommand{\floor}[1]{\lfloor{#1}\rfloor}

\renewcommand{\Pr}{\mathbf{P}\,}

\DeclareMathOperator{\Li}{Li}

\newtheorem{theorem}{Theorem}
\newtheorem{proposition}[theorem]{Proposition}
\theoremstyle{remark}
\newtheorem{remark}[theorem]{Remark}
\newtheorem{example}[theorem]{Example}
\numberwithin{theorem}{section}

\numberwithin{equation}{section}

\title[Uniform $m$-ary search trees]{A repertoire for
  additive functionals of uniformly distributed  $m$-ary~search trees}
\author{James Allen Fill}
\author{Nevin Kapur}

\thanks{The research of James Allen Fill was supported by \textsc{NSF}
  Grants \textsc{DMS}--0104167 and \textsc{DMS}--0406104, and by The
  Johns Hopkins University's Acheson J.~Duncan Fund for the
  Advancement of Research in Statistics.  Nevin Kapur's research was
  partially supported by \textsc{NSF} Grant~0049092, and the by Center
  for Mathematics of Mathematics of Information at the California
  Institute of Technology.}

\address[James Allen Fill]{Department of Applied Mathematics and
  Statistics, The 
  Johns Hopkins University, 3400 N.~Charles St., Baltimore MD
  21218-2682}
\urladdr{http://www.ams.jhu.edu/\~{}fill/}
\email{jimfill@jhu.edu}

\address[Nevin Kapur]{Department of Computer Science, California Institute of
  Technology, MC 256-80, 1200 E.~California Blvd., Pasadena CA 91125}
\urladdr{http://www.cs.caltech.edu/\~{}nkapur/}
\email{nkapur@cs.caltech.edu}



\date{February~24, 2005}

\begin{document}
\maketitle

\begin{abstract}
  Using recent results on singularity analysis for Hadamard products
  of generating functions, we obtain the limiting distributions
  for additive functionals on $m$-ary search trees on $n$~keys with
  toll sequence (i)~$n^\alpha$ with $\alpha \geq 0$ ($\alpha=0$ and
  $\alpha=1$ correspond roughly to the space requirement and
  total path length, respectively); (ii)~$\ln
  \binom{n}{m-1}$, which corresponds to the so-called shape
  functional; and (iii)~$\mathbf{1}_{n=m-1}$, which corresponds to the
  number of leaves.
\end{abstract}

\section{Introduction}
\label{sec:introduction}

We begin by providing a brief overview of $m$-ary search trees.
For integer $m \geq 2$, the $m$-ary 
search tree, or multiway tree, generalizes the binary search tree.
The quantity~$m$ is called the \emph{branching factor}.
According to~\cite{MR93f:68045}, search trees with branching factors
higher than~2 were first suggested by Muntz and
Uzgalis~\cite{muntz71:_dynam} ``to solve internal memory problems with
large quantities of data.''  For further background we refer the reader
to~ \cite{knuth97,knuth98} and~\cite{MR93f:68045}.

We consider the space of \emph{$m$-ary search trees} on $n$~keys, and
assume that the keys can be linearly ordered.  Since we shall be
concerned only with the structure of the tree and not its specific
contents, we can then without loss of generality take the set of keys
to be $[n] := \{1,2,\ldots,n\}$. An $m$-ary search tree can be
constructed from a sequence $s$ of $n$ distinct keys in the following
way:
\begin{enumerate}[(a)]
\item If $n < m$, then all the keys are stored in the root node in
  increasing order. \label{item:m-ary-def1}
\item If $n \geq m$, then the first $m-1$ keys in the sequence are
  stored in the root in increasing order, and the remaining $n-(m-1)$
  keys are stored in the $m$~subtrees subject to the condition that if
  $\kappa_1 < \kappa_2 < \cdots < \kappa_{m-1}$ denotes the ordered
  sequence of keys in the root, then the keys in the $j$th subtree are
  those that lie between $\kappa_{j-1}$ and $\kappa_{j}$, where
  $\kappa_0 := 0$ and $\kappa_{m} := n+1$, sequenced as in~$s$.
  \label{item:m-ary-def2}
\item Recursively, all the subtrees are $m$-ary search trees that satisfy
  conditions~(\ref{item:m-ary-def1}),~(\ref{item:m-ary-def2}),
  and~(\ref{item:m-ary-def3}). \label{item:m-ary-def3}
\end{enumerate}
In this work we consider additive functionals on $m$-ary search trees,
as we describe next.

Fix $m \geq 2$.  Given an $m$-ary search tree~$T$, let $L_1(T), \dots,
L_m(T)$ denote the subtrees rooted at the children of the root
of~$T$.  The \emph{size}~$|T|$ of a tree~$T$ is the number of keys in it.
We will call a  functional~$f$ on~\( m \)-ary search trees
\emph{additive} if it satisfies the recurrence
\begin{equation}
  \label{eq:2.1}
  f(T) = \sum_{i=1}^m f( L_i(T) ) + b_{|T|},
\end{equation}
for any tree~\( T \) with~\( |T| \geq m -1 \). Here~\( (b_n)_{n \geq
  m-1} \) is a given sequence, henceforth called the \emph{toll
  sequence} or \emph{toll function}.
Note that the recurrence~\eqref{eq:2.1} does not make any reference
to~\( b_n \) for~\( 0 \leq n \leq m-2 \) nor specify the initial
  conditions \( f(T) \) for~\(
0 \leq |T| \leq m-2 \).

Several interesting examples can be cast as additive functionals.
\begin{example}
  \label{example:space-requirement} If we specify~\( f(T) \)
  arbitrarily for~\( 0 \leq |T| \leq m-2 \) and take~\( b_n \equiv c
  \) for~\( n \geq m-1 \), we obtain the ``additive functional''
  framework of~\cite[\S3.1]{MR93f:68045}.  (Our definition of an
  additive functional substantially generalizes this notion.)  In
  particular if we define~\( f(\emptyset) := 0 \) and~\( f(T) := 1 \)
  for the unique \( m \)-ary search tree~\( T \) on~\( n \) keys for
  \( 1 \leq n \leq m-2 \) and let \( b_n \equiv 1 \) for \( n \geq m-1
  \), then \( f(T) \) counts the number of nodes in~$T$ and thus gives
  the \emph{space requirement} functional discussed
  in~\cite[\S3.4]{MR93f:68045}.
\end{example}
\begin{example}
  \label{example:leaves}
  If we define $f(T) := 0$ when $|T|=0$, $f(T) := 1$ when $1 \leq |T|
  \leq m-2$, and $b_n := \mathbf{1}_{n=m-1}$, then $f$ is the
  \emph{number of leaves} in the $m$-ary search tree.
\end{example}
\begin{example}
  \label{example:path-length}
  If we define \( f(T) := 0 \) when \( 0 \leq |T| \leq m-2 \) and
  \( b_n := n-(m-1) \) for \( n \geq m-1 \), then \( f \) is the
  \emph{internal path length} functional discussed
  in~\cite[\S3.5]{MR93f:68045}:\ $f(T)$ is the sum of all root-to-key
  distances in $T$.
  
  In this work we choose to treat explicitly the toll~$n$, rather
  than~$n-(m-1)$.  However our techniques reveal that the lead-order
  asymptotics of moments and the limiting distributions of these two
  additive functionals are the same.
\end{example}
\begin{example}
  \label{example:shape-functional}
  As described above, each permutation of $[n]$ gives rise to an
  $m$-ary search tree.  Suppose we place the uniform distribution on
  such permutations. This induces a distribution on $m$-ary search
  trees called the \emph{random permutation model}. Denote its
  probability mass function by~$Q$.  
  Dobrow and Fill~\cite{MR97k:68038} noted  that
  \begin{equation}
    \label{eq:2.2}
    Q(T) = \frac{1}{\prod_{x} \binom{|T_x|}{m-1}},
  \end{equation} where the product in~\eqref{eq:2.2} is over all
  nodes in~$T$ that contain $m-1$ keys.  This functional is sometimes
  called the ``shape functional'' as it serves as a crude measure of
  the ``shape'' of the tree, with ``full'' trees (such as the complete
  tree) achieving the larger values of~$Q$.  For further discussions
  along these lines, consult~\cite{MR97k:68038}
  and~\cite{MR97f:68021}.  If we define \( f(T) := 0\) for \( 0 \leq
  |T| \leq m-2 \) and \( b_n := \ln{\binom{n}{m-1}} \) for \( n \geq
  m-1 \), then \( f(T) = -\ln{Q(T)} \).  Henceforth throughout this
  paper we will refer to $-\ln{Q}$ (rather than
  $Q$) as the \emph{shape functional}.  
\end{example}

Several authors~\cite{MR90a:68012,MR95i:68030,MR1871558,FK-transfer}
have studied additive functionals under the random permutation model.
Clearly the random permutation model does not induce the uniform
distribution on $m$-ary search trees with $n$~keys since different
permutations can give rise to the same tree.  In this paper we
consider additive functionals under the uniform model, i.e., when each
tree on $n$~keys is considered equally likely.  The shape functional
for the case $m=2$ (uniformly distributed binary search trees) was
considered by Fill~\cite{MR97f:68021}, who derived (limited)
asymptotic information about its mean and variance.  Limiting
distributions for the shape functional and other additive functionals
treated in the present paper were identified
in~\cite{fill03:_limit_catal} for $m=2$.  We now generalize these
results to include all values of~$m$.  What makes the analysis for
general~$m$ significantly more intricate is that several key
quantities (such as the number~$\rho$ discussed at the beginning of
Section~\ref{sec:singular-expansions}) are for general~$m$ known only
implicitly.

One motivation for the present paper can be understood in the context
of the shape functional.  The probability mass function~$Q$
corresponding to the random permutation model (a reasonably realistic
model in practice) is an object of natural interest.  Dobrow and
Fill~\cite{MR97k:68038} determined the smallest and largest values
of~$Q$; but what are ``typical'' values?  We can study this question
probabilistically by placing a distribution on~$T$ and considering the
distribution of~$Q(T)$.  Two rather natural choices for this
distribution are $Q$ itself (as treated in~\cite{FK-transfer}) and the
uniform distribution on trees (as treated herein).


We follow the ``repertoire'' approach of Greene and
Knuth~\cite{MR92c:68067}, determining the effect of a family of basic
tolls (for example, those of the form $n^\alpha$).  Then the effect of
a new toll could be determined by expressing it in terms of the basic
tolls.

For tolls of the form $n^\alpha$ with $\alpha \geq 0$ and the tolls
$\ln{\binom{n}{m-1}}$ and $\mathbf{1}_{n=m-1}$, we determine
asymptotics of moments of all orders and our main results
(Theorems~\ref{thm:dist-ne-12}, \ref{thm:alpha=half},
\ref{thm:shape-clt}, \ref{thm:space-clt}, and~\ref{thm:leaves-clt})
use these to yield limiting distributions.  Here, in broad terms for
the toll~$n^\alpha$, is a summary of lead-order results under both the
random permutation model and the uniform model:
\begin{table}[htbp]
  \centering
  \begin{tabular}{c|cc}
    & \multicolumn{2}{c}{Model} \\
    Toll function~$n^\alpha$ & Random permutation & Uniform \\
    \hline
    $\alpha$ smaller than~1/2 & $n$ & $n$ \\
    $\alpha$ between 1/2 and 1 & $n$ &  $n^{\alpha+\frac12}$ \\
    $\alpha$ bigger than 1 & $n^\alpha$ & $n^{\alpha+\frac12}$
  \end{tabular}
  \caption{Order of magnitude of the additive functional
    corresponding to the toll~$n^\alpha$.}
  \label{tab:typical}
\end{table}\\
It is not surprising that the orders of magnitude under the uniform
model are at least as large as under the random permutation model.
Indeed, it is well known that trees produced by the uniform model are
generally much ``stringier'' than trees produced by the random
permutation model; for example,  height is of order~$\sqrt n$ under
the uniform model and order~$\log n$ under the random permutation
model.  Furthermore ``stringy'' trees tend to give large values of
the functional.

Qualitatively the uniform model differs significantly from the random
permutation model, where, for example, there is a ``phase change'' in
the limiting behavior at $m=26$ from asymptotic normality to
non-existence of a limiting distribution, for any toll whose order of
growth does not exceed~$n^{1/2}$; see~\cite{FK-transfer} for precise
results.  On the other hand, for all~$m$ the uniform model leads to
the normal distribution for the shape functional, space requirement,
and number of leaves, and to (apparently) non-normal distributions for
tolls of the form~$n^\alpha$ with $\alpha > 0$.

We use methods from analytic combinatorics, in particular singularity
analysis of generating functions~\cite{MR90m:05012}, to derive the
asymptotics of moments of the functional under consideration and then
the method of moments to characterize the limiting distribution.  A
key singularity analysis tool is the newly-developed ``Zigzag
algorithm''~\cite{FFK} to handle Hadamard products of generating
functions.

The limiting distributions (and even local limit theorems) for the
space requirement and the number of leaves presumably can also be
derived using Theorem~2 of~\cite{MR1284406} since the bivariate
generating function for these parameters satisfy suitable functional
equations.  (This is \emph{not} the case for the other tolls that we
consider.)  We include our proofs of these results for completeness
and uniformity of treatment of tolls.

The paper is organized as follows.  In Section~\ref{sec:preliminaries}
we set up the problem using generating functions.  In
Section~\ref{sec:singular-expansions}, a singular expansion for the
generating function of the number of $m$-ary search trees on $n$~keys
is obtained.  Sections~\ref{sec:toll-nalpha},
\ref{sec:shape-functional}, \ref{sec:space-requirement},
and~\ref{sec:number-leaves} derive limiting distributions for the
additive functionals corresponding to the tolls $n^{\alpha}$ ($\alpha
> 0$), $\ln{\binom{n}{m-1}}$ (shape functional), $1$ (space
requirement), and $\mathbf{1}_{n=m-1}$ (number of leaves),
respectively.

\emph{Notation}.  Throughout, we will use $[z^n]f(z)$ to denote the
coefficient of $z^n$ in the Taylor series expansion of $f(z)$ around
$z=0$.  We use $\mathcal{L}(Y)$ to denote the law (or distribution) of
a random variable~$Y$, the symbol $\eqinlaw$ to denote equality in
law, and $\toinlaw$ to denote convergence in law.  We denote the
(univariate) normal distribution with mean~$\mu$ and
variance~$\sigma^2$ by~$N(\mu, \sigma^2)$.

\section{Preliminaries}
\label{sec:preliminaries}
Our starting point is the recursive construction of $m$-ary search
trees.  Let $X_n \equiv X_n(T)$ denote an additive functional on a
random $m$-ary search tree~$T$ on $n$~keys. Let $\mathbf{J} \equiv
(J_1, \dots, J_m$) be the (random) vector of sizes of the subtrees
rooted at the children of the root of~$T$.  If $T$ is a uniformly
distributed $m$-ary search tree on $n$~keys, then $X_n$ satisfies the
distributional recurrence
\begin{equation}
  \label{eq:1.1}
  X_n \eqinlaw \sum_{k=1}^m X_{J_k}^{(k)} + b_n, \quad n \geq m-1,
\end{equation}
with $(X_0,\dots,X_{m-2}) =: \mathbf{x}$ denoting the vector of
deterministic values of the functional for trees with fewer than $m-1$
keys.  The sequence $(b_n)_{n \geq m-1}$ is called the \emph{toll
  sequence}.  In~\eqref{eq:1.1}, $\eqinlaw$ denotes equality in law (i.e.,\ in
distribution), and on the right,
\begin{itemize}
\item for each $k = 1, \dots, m$, we have $X_j^{(k)} \eqinlaw X_j$;
\item the quantities $\mathbf{J}$; $X_0^{(1)}, \dots, X_{n-(m-1)}^{(1)}$;
  $X_0^{(2)}, \dots, X_{n-(m-1)}^{(2)}$; \ldots ; $X_0^{(m)}, \dots,
  X_{n-(m-1)}^{(m)}$ are all independent; \sloppy
\item the distribution of~$\mathbf{J}$ if given by
\begin{equation}
  \label{eq:55}
  \Pr[ J_1=j_1, \ldots, J_m=j_m ] = \frac{\tau_{j_1} \cdots
  \tau_{j_m}}{\tau_n},
\end{equation}
for $(j_1,\dots,j_m) \geq \mathbf{0}$ with $j_1 + \dots + j_m = n-(m-1)$,
where $\tau_k \equiv \tau_k(m)$ is the number of $m$-ary search trees
on $k$ keys.
\end{itemize}
(Throughout we will take $m \geq 2$ to be fixed and so will suppress the
dependence of various parameters on~$m$.)

Denote the $s$th moment of $X_n$ by $\mun{s} := \E X_n^s$.  Now taking
the $s$th power of~\eqref{eq:1.1} and conditioning on $(J_1,\dots,J_m)$
gives
\begin{equation*}
  \mun{s} = \sum_{s_0 + \dots + s_m = s} \binom{s}{s_0,\dots,s_m}
  b_n^{s_0} \sideset{}{^*}\sum \frac{\tau_{j_1} \cdots
  \tau_{j_m}}{\tau_n} \mu_{j_1}^{[s_1]} \cdots \mu_{j_m}^{[s_m]},
\end{equation*}
where $\sideset{}{^*}\sum$ denotes the sum over all $m$-tuples
$(j_1,\dots, j_m) \geq \mathbf{0}$ such that $\sum_{i=1}^m j_i
= n - (m-1)$.  Isolating the terms in the sum where $s_i=s$ for some~$i
\in [m]$, we get
\begin{align}
  \tau_n \mun{s} &= m \sideset{}{^*} \sum \tau_{j_1} \cdots \tau_{j_m}
  \, \mu_{j_1}^{[s]} \notag \\
  & \quad {}+ \sum_{\substack{s_0 + \dots + s_m = s \\ s_1,
  \dots, s_m < s}} \binom{s}{s_0, \dots, s_m} b_n^{s_0}
  \sideset{}{^*}\sum \tau_{j_1}\mu_{j_1}^{[s_1]} \cdots
  \tau_{j_m} \mu_{j_m}^{[s_m]} \notag                       \\
  \label{eq:2.4}
                 &= m \sum_{j_1=0}^{n-(m-1)} \tau_{j_1}
  \mu_{j_1}^{[s_1]} \sum_{ j_2 + \dots + j_m = n - (m-1)-j_1} \tau_{j_2}
  \cdots \tau_{j_m} + r_n^{[s]},
\end{align}
where
\begin{equation}\label{eq:47}
  r_n^{[s]} := \sum_{\substack{s_0 + \dots + s_m = s\\ s_1,\dots,s_m < s}}
  \binom{s}{s_0, \dots, s_m} b_n^{s_0} \sideset{}{^*}\sum \tau_{j_1}
  \mu_{j_1}^{[s_1]} \cdots \tau_{j_m} \mu_{j_m}^{[s_m]}.
\end{equation}
Introduce generating functions
\begin{equation*}
  \mu^{[s]}(z) := \sum_{n=0}^\infty \tau_n \mu_n^{[s]} z^n, \quad
  r^{[s]}(z) := \sum_{n=0}^\infty r_n^{[s]} z^n, \quad \tau(z) :=
  \sum_{n=0}^\infty \tau_n z^n.
\end{equation*}
Multiplying~\eqref{eq:2.4} by $z^n$ and summing over $n \geq m-1$
yields (observe that $\tau_0 = \dots = \tau_{m-2} = 1$ and $r_0^{[s]}
= \dots = r_{m-2}^{[s]} = 0$)
\begin{equation*}
  \mu^{[s]}(z) - \sum_{j=0}^{m-2} x_j^s z^j = m z^{m-1} \mu^{[s]}(z)
  \tau^{m-1}(z) + r^{[s]}(z),
\end{equation*}
so that
\begin{equation}
  \label{eq:3.7}
  \mu^{[s]}(z) = \frac{r^{[s]}(z) + \sum_{j=0}^{m-2} x_j^s z^j}{1 -
  m[z \tau(z)]^{m-1}}.
\end{equation}
Furthermore
\begin{align}
  r^{[s]}(z) &= \sum_{\substack{s_0 + \dots + s_m = s\\ s_1, \dots,
  s_m < s}} \binom{s}{s_0, \dots, s_m} \sum_{n=0}^\infty b_n^{s_0}
  \left\{[z^{n-(m-1)}] \left[\muz{s_1} \cdots \muz{s_m}\right] \right\}
  z^n \notag \\
  \label{eq:2A.4}
  &= \sum_{\substack{s_0 + \dots + s_m = s\\ s_1, \dots,
      s_m < s}} \binom{s}{s_0, \dots, s_m} b^{\odot s_0} (z) \odot
  \left( z^{m-1} \muz{s_1} \cdots \muz{s_m} \right),
\end{align}
where
\begin{equation*}
  b(z) := \sum_{n=0}^\infty b_n z^n
\end{equation*}
and $f(z) \odot g(z) \equiv (f \odot g)(z)$ is the Hadamard product of
the power series~$f$ and~$g$.  Note that since $[z^n] \left( z^{m-1}
  \muz{s_1} \cdots \muz{s_m} \right) = 0$ for $0 \leq n \leq m-2$ we
may instead use
\begin{equation*}
  b(z) := \sum_{n=m-1}^\infty b_n z^n
\end{equation*}
when convenient.

\section{Singular expansions}
\label{sec:singular-expansions}

We will employ singularity
analysis~\cite{MR90m:05012,MR2000a:05015,FFK} to derive asymptotics of
$\mu_n^{[s]}$ using~\eqref{eq:3.7}.  In order to do so we need a
singular expansion for $\tau(z)$ around its dominant singularity.  We
will use the theory of analytic continuation of algebraic functions
(see, for example,~\cite[\S III.45]{markushevich}
or~\cite[\S{}VII.4]{flajolet-sedgewick}) to derive such an
expansion.  The terminology used is
from~\cite[\S{}VII.4]{flajolet-sedgewick}.

Before we begin, we note that Fill and Dobrow~\cite{MR99f:68034} were
able to use large-deviations techniques to obtain lead-order
asymptotics of~$\tau_n$.  However their techniques do not seem to be
sufficient to derive the higher-order results we will need.

We now proceed with our analytic approach.  As observed by Fill and
Dobrow~\cite{MR99f:68034}, it follows from the recursive definition of
$m$-ary search trees that
\begin{equation}
  \label{eq:3}
  \tau(z) - \sum_{j=0}^{m-2} z^j = z^{m-1} \tau^m(z).
\end{equation}
Thus $\tau(z)$ is an algebraic series satisfying $P(z, \tau(z)) = 0$, where
\begin{equation}
  \label{eq:48}
  P(z,w) := z^{m-1} w^m - w + \sum_{j=0}^{m-2} z^j.
\end{equation}
The exceptional set of~$P$ [excluding $z=0$, at which  $\tau(z)$
clearly has no singularity] is 
\begin{align*}
  & \bigcup_{w \in \mathbb{C}} \Bigl\{ z\!: P(z,w) = 0 \text{ and }
  \frac{\partial}{\partial w} P(z,w) = 0 
  \Bigr\} \\
  & = \bigcup_{w \in \mathbb{C}} \Bigl\{ z\!: z^{m-1} w^m - w +
  \sum_{j=0}^{m-2} z^j = 0
  \text{ and } m(zw)^{m-1} - 1 = 0 \Bigr\} \\
  &= \left\{ z\!: m^m \left( \sum_{j=1}^{m-1} z^j
  \right)^{m-1} = (m-1)^{m-1}   \right\}.
\end{align*}
The singularities of~$\tau(z)$ lie in the exceptional set.
It is clear~\cite[Theorem~3.1]{MR99f:68034} that there exists a
unique $\rho \in (0,1)$ contained in this set.  Furthermore, since the Taylor
coefficients of $\tau(z)$ are nonnegative, by Pringsheim's
theorem~\cite[Theorem~I.17.13]{markushevich}, $\rho$ is a dominant
singularity of~$\tau(z)$.
It is straightforward to check that the polynomial system given by
writing $P(z,w) = 0$ in the form $w = \Phi(z,w)$ is a-proper,
a-positive, a-irreducible, and a-aperiodic (cf.~\cite[\S
VII.4.2]{flajolet-sedgewick}), so that by Theorem~VII.7
of~\cite{flajolet-sedgewick} we have that $\rho$ is the unique
dominant singularity and as $z \to \rho$ a singular expansion of the
form
\begin{equation}
  \label{eq:1}
  \tau(z) \sim \sum_{l \geq 0} a_l (1-\rho^{-1}z)^{{l}/2}.
\end{equation}
\begin{remark}
  \label{rem:tau_n_asymptotics}
  Singularity analysis immediately yields from~\eqref{eq:1} a complete
  asymptotic expansion for $\tau_n$, the number of $m$-ary search
  trees on $n$~keys:
  \begin{equation}
    \label{eq:13}
    \tau_n \sim \rho^{-n} \sum_{l \geq 0} \frac{a_{2l+1}}{\Gamma(-l-\frac12)}
    n^{-l-\frac32}.
  \end{equation}
  In particular,
  \begin{equation*}
    \tau_n = [1 + O(n^{-1})] \frac{-a_1}{2\sqrt\pi} n^{-3/2}
    \rho^{-n}.
  \end{equation*}
\end{remark}

\subsection{Determination of the coefficients $a_l$}
\label{sec:determination-a_ls}

Define $w_\rho := \frac{m}{m-1} \sum_{j=0}^{m-2} \rho^j$, so that
\begin{equation*}
  P(\rho, w_\rho) = 0 \text{ and } \frac{\partial}{\partial w} P(\rho,
  w) \biggr\vert_{w=w_\rho} = 0.
\end{equation*}
Using the definition of~$\rho$ and the fact that $w_\rho > 0$ by
definition, we have $w_\rho = m^{-\frac{1}{m-1}} \rho^{-1}$.  Now
$\frac{\partial}{\partial w}P(\rho,w)$ is negative, zero, or positive
as $w > 0$ is less than, equal to, or greater than~$w_\rho$.  Hence,
for $w > 0$, $P(\rho, w) = 0$ if and only if $w=w_\rho$.  But $a_0 >
0$ and $0 =
P(\rho, \tau(\rho)) = P(\rho, a_0)$, so that
\begin{equation}
  \label{eq:4}
  a_0 = w_\rho =  m^{-\frac1{m-1}} \rho^{-1}.
\end{equation}

To obtain values of $a_l$ for $l \geq 1$, we rewrite~\eqref{eq:3}
for $z \ne 1$ as
\begin{equation*}
  z^{m-1} \tau^m(z) - \tau(z) + \frac{1-z^{m-1}}{1-z} = 0
\end{equation*}
and, then defining $Z := 1-\rho^{-1}z$, equivalently as
\begin{equation}
  \label{eq:56}
  1 + \rho^{m-1} (1-Z)^{m-1} \bigl[(1-\rho+\rho Z)\tau^m(z) - 1 \bigr]
  - (1-\rho + \rho Z) \tau(z).
\end{equation}
By comparing the coefficients of~$Z$ in this equation and observing
that $a_1 < 0$ we obtain
\begin{equation}
  \label{eq:50}
  a_1 = -\sqrt{2m\alpha^*} m^{-\frac{m}{m-1}} \rho^{-1}, 
\end{equation}
where, matching the notation of~\cite{MR99f:68034}, we define the key
quantity
\begin{equation}
  \label{eq:18}
  \alpha^* := m - \left(m^\frac{m}{m-1} - 1\right)(\rho^{-1} -
  1)^{-1}.
\end{equation}
In the sequel we will also need the following relation, which follows
from comparing coefficients of $Z^{3/2}$ in~\eqref{eq:56}:
\begin{equation}
  \label{eq:49}
  \frac{a_0(a_0 - a_2)}{a_1^2} = \frac{m-2}6.
\end{equation}

Let $\mathcal{A}$ denote a generic (formal) power series in~$Z$,
possibly different at each appearance.  Similarly, let $\mathcal{P}_d$
denote a generic polynomial in~$Z$ of degree at most~$d$.  In the
sequel we will likewise use $\mathcal{N}$ to denote a generic (formal)
power series in powers of $n^{-1}$.  Then, using~\eqref{eq:1}
and~\eqref{eq:4}, we have
\begin{equation}
  \label{eq:5}
    (1 - m [z\tau(z)]^{m-1})^{-1} \sim \frac{a_0}{-a_1(m-1)}Z^{-1/2} +
  c_0 + Z^{1/2}\A + Z\A,
\end{equation}
where, using~\eqref{eq:49}, we have 
\begin{equation}
  \label{eq:57}
  c_0 := \frac{m-2}{3(m-1)};
\end{equation}
\begin{equation}
  \label{eq:11}
  z^{m-1} \tau^m(z) \sim a_0 m^{-1} + a_1 Z^{1/2} + Z\mathcal{A} +
  Z^{3/2}\mathcal{A};
\end{equation}
and
\begin{equation}
  \label{eq:6}
  \sum_{j=0}^{m-2} x_j^s z^j = \sum_{j=0}^{m-2} x_j^s \rho^j + Z
  \mathcal{P}_{m-3}.
\end{equation}
Thus, by singularity analysis,
\begin{equation}
  \label{eq:14}
  [z^n] \left[ z^{m-1} \tau^m(z) \right] \sim n^{-3/2} \rho^{-n}
  \left( \frac{-a_1}{2\sqrt\pi} + n^{-1} \mathcal{N} \right).
\end{equation}

\subsection{Generalized polylogarithms}
\label{sec:gener-polyl}

  For $\alpha$ an arbitrary complex number and $r$ a nonnegative
  integer, the \emph{generalized polylogarithm} function
  $\Li_{\alpha,r}$ is defined  for $|z| < 1$ by
  \begin{equation*}
    \Li_{\alpha,r}(z) := \sum_{n=1}^\infty \frac{(\ln{n})^r}{n^\alpha} z^n.
  \end{equation*}
  We record here three singular expansions (as $z \to \rho$) that are
  computed using singularity analysis of generalized
  polylogarithms~\cite{MR2000a:05015}:
\begin{equation}
  \label{eq:15}
  \begin{split}
    \Li_{3/2,0}(\rho^{-1}z) &\sim \zeta(3/2) - 2\sqrt\pi Z^{1/2} +
    Z\A + Z^{3/2}\A, \\
    \Li_{3/2,1}(\rho^{-1}z) &\sim -\zeta'(3/2) - 2\sqrt\pi Z^{1/2}
    \ln{Z^{-1}} - 2\sqrt\pi[ 2(1 - \ln 2) - \gamma] Z^{1/2} \\
    & \quad {}+ Z\A + (Z^{3/2} \log Z)\A + Z^{3/2}\A, \\
    \Li_{1,1}(\rho^{-1}z) & \sim \frac12 \ln^2 Z^{-1} - \gamma \ln
    Z^{-1} + \A + (Z\log Z)\A.
  \end{split}
\end{equation}
These expansions will be utilized in the analysis of the shape functional in
Section~\ref{sec:shape-functional}.

\subsection{Zigzag algorithm}
\label{sec:zigzag-algorithm}
For the reader's convenience we present the Zigzag algorithm, which is
used extensively in the rest of this paper to determine singular
expansions of Hadamard products.  The validity of the algorithm was
established recently in~\cite{FFK}, to which the reader is referred
for further background discussion.
\begin{itemize}\def\Asympt{\operatorname{Asympt}}
\item[] \emph{\bf ``Zigzag'' Algorithm}. [Computes the singular
  expansion of
  $f\odot g$ up to $O(|1-z|^C)$. ]\\[1truemm]
  {\bf 1.} Use singularity analysis to determine separately the
  asymptotic expansions of~$f_n=[z^n]f(z)$ and~$g_n=[z^n]g(z)$ into
  descending powers of~$n$.
  \\
  {\bf 2.} Multiply the resulting expansions and reorganize to obtain
  an asymptotic expansion for the product~$f_n g_n$.
  \\
  {\bf 3.} Choose a basis $\mathcal{B}$ of singular functions, for
  instance, the standard basis $\mathcal{B}=\left\{(1-z)^\beta
    [\ln(1-z)]^k\right\}$, or the polylogarithm basis $\mathcal{B}=\left\{
    \Li_{\beta,k}(z)\right\}$.  Construct a function~$H(z)$ expressed
  in terms of~$\mathcal{B}$ whose singular behavior is such that the
  asymptotic form of its coefficients~$h_n$ is compatible with that
  of~$f_n g_n$ up to the needed error terms.
  \\
  {\bf 4.} Output the singular expansion of $f\odot g$ as the quantity
  $H(z)+P(z)+O(|1-z|^C)$, where $P$ is a polynomial in~$(1-z)$ of
  degree less than~$C$.
\end{itemize}
The reason for the addition of a polynomial in Step~{\bf 4}
is that integral powers of~$(1-z)$ do not leave a trace
in coefficient asymptotics since their contribution is 
asymptotically null. 
The Zigzag Algorithm is principally useful for determining the
divergent part of expansions.  If needed, the coefficients in the
polynomial~$P$ can be expressed as values of the function $f\odot g$
and its derivatives at~1 once it has been stripped of its
nondifferentiable terms.

\section{The toll~$n^\alpha$}
\label{sec:toll-nalpha}

The main theorems of this section are Theorems~\ref{thm:dist-ne-12}
and~\ref{thm:alpha=half}, which give limiting distributions for the
additive functionals corresponding to the tolls $n^{\alpha}$ with
$\alpha >0$. Although the normalization required to produce a limiting
distribution depends on~$m$, our results exhibit a striking
\emph{invariance principle}: The limiting distributions themselves do
not depend on the value of~$m$ (and thus in particular, have already
arisen when $m=2$ in~\cite{fill03:_limit_catal}).

\subsection{Mean}
\label{sec:mean}

We consider the mean for the toll $b_n \equiv n^\alpha$, where $\alpha >
0$.  Using $s=1$ in~\eqref{eq:2A.4} we have
\begin{equation*}
  r^{[1]}(z) = b(z) \odot \left[ z^{m-1} \tau^m(z) \right],
\end{equation*}
and consequently, by~\eqref{eq:47} and~\eqref{eq:13},
\begin{equation}\label{eq:40}
  [z^n] r^{[1]}(z) = r_n^{[1]} = b_n \tau_n \sim 
  n^{\alpha - \frac32} \rho^{-n} \left( \frac{-a_1}{2\sqrt\pi} +
  n^{-1} \mathcal{N} \right).
\end{equation}

Until further notice, assume that $\alpha \not\in
\{1/2,3/2,\dots\}$.  (The contrary cases are considered later
in this section.)  We employ the Zigzag Algorithm outlined in
Section~\ref{sec:zigzag-algorithm}. A compatible singular expansion
for $r^{[1]}(z)$ is given by
\begin{equation}
  \label{eq:7}
  r^{[1]}(z) \sim \frac{-a_1}{2\sqrt\pi} \Gamma(\alpha - \tfrac12)
  Z^{-\alpha + \frac12} + Z^{-\alpha + \frac32}\mathcal{A} + \mathcal{A}.
\end{equation}
If $\alpha > 1/2$, then using~\eqref{eq:5},~\eqref{eq:6},
and~\eqref{eq:7} in~\eqref{eq:3.7} we obtain
\begin{equation}
  \label{eq:10}
  \muz{1} \sim \frac{a_0 \Gamma(\alpha - \frac12)}{2\sqrt\pi(m-1)}
  Z^{-\alpha} + Z^{-\alpha + \frac12}\mathcal{A} + Z^{-\alpha +
  1}\mathcal{A} + Z^{-1/2}\mathcal{A} + \mathcal{A},
\end{equation}
whence, by singularity analysis,
\begin{equation*}
  \rho^{n} \mun{1} \tau_n \sim \frac{a_0
  \Gamma(\alpha-\frac12)}{2\sqrt\pi(m-1)\Gamma(\alpha)} n^{\alpha-1} +
  n^{\alpha-\frac32}\mathcal{N} + n^{\alpha-2}\mathcal{N} +
  n^{-1/2}\mathcal{N}.
\end{equation*}
The singular expansion for $\tau_n$ at~\eqref{eq:13} then gives
\begin{equation*}
  \mu_n^{[1]} \sim 
  \frac{a_0\Gamma(\alpha-\frac12)}{(-a_1)(m-1)\Gamma(\alpha)} 
  n^{\alpha + \frac12} +
  n^{\alpha} \mathcal{N} + n^{\alpha-\frac12}\mathcal{N} + n\mathcal{N}.
\end{equation*}
On the other hand, if $\alpha < 1/2$, the dominant term
in~\eqref{eq:7} is now the constant term so that
\begin{equation}
  \label{eq:46}
  r^{[1]}(z) + \sum_{j=0}^{m-2} x_j z^j \sim C_\alpha +
  \frac{-a_1}{2\sqrt\pi} \Gamma(\alpha - \tfrac12)
  Z^{-\alpha + \frac12}  + Z\mathcal{A} + Z^{-\alpha +
  \frac32}\mathcal{A},
\end{equation}
where
\begin{equation}
  \label{eq:41}
  C_\alpha := r^{[1]}(\rho) + \sum_{j=0}^{m-2} x_j \rho^j  =
  \sum_{n=m-1}^\infty \rho^n n^\alpha \tau_n + \sum_{j=0}^{m-2} x_j
  \rho^j.
\end{equation}
Then, using~\eqref{eq:5} and~\eqref{eq:46}
in~\eqref{eq:3.7} we obtain
\begin{equation}
  \label{eq:52}
  \muz{1} \sim \frac{a_0 C_\alpha}{(m-1)(-a_1)} Z^{-1/2} + \frac{a_0
  \Gamma(\alpha - \frac12)}{2\sqrt\pi(m-1)} 
  Z^{-\alpha} + \mathcal{A} + Z^{-\alpha + \frac12}\mathcal{A} +
  Z^{-\alpha + 1}\mathcal{A} + Z^{1/2}\mathcal{A},
\end{equation}
whence singularity analysis and the singular expansion of $\tau_n$
yields
\begin{equation}  \label{eq:9}
  \mun{1} \sim \frac{2 a_0 C_\alpha}{(m-1)a_1^2} n +
  \frac{a_0\Gamma(\alpha-\frac12)}{(-a_1)(m-1)\Gamma(\alpha)}  
  n^{\alpha + \frac12} +
  n^{\alpha} \mathcal{N}  + \mathcal{N} +
  n^{\alpha-\frac12}\mathcal{N}.
\end{equation}

If $\alpha \in \{3/2, 5/2, \dots \}$, then logarithmic terms
appear in the singular expansion compatible with~\eqref{eq:40}, so that
\begin{equation*}
  r^{[1]}(z) \sim \frac{-a_1}{2\sqrt\pi} \Gamma(\alpha-\tfrac12)
  Z^{-\alpha + \frac12} + Z^{-\alpha + \frac32}\A + (\log Z)\A.
\end{equation*}
This leads to
\begin{equation}
  \label{eq:39}
  \muz{1} \sim \frac{a_0 \Gamma(\alpha-\frac12)}{2\sqrt\pi(m-1)}
  Z^{-\alpha} + Z^{-\alpha+\frac12}\A + Z^{-\alpha+1}\A +
  (Z^{-1/2}\log Z)\A + (\log Z)\A
\end{equation}
and consequently
\begin{equation*}
  \mun{1} \sim \frac{a_0
  \Gamma(\alpha-\frac12)}{-a_1(m-1)\Gamma(\alpha)} n^{\alpha+\frac12}
  + n^{\alpha}\N + n^{\alpha-\frac12}\N + (n \log n)\N.
\end{equation*}
Observe that the lead-order term and the order of growth of the
remainder [$O(|Z|^{-\alpha+\frac12})$]  in
the expansion of~$\muz{1}$ at~\eqref{eq:39} are the same as
at~\eqref{eq:10}.

Finally, we consider $\alpha=1/2$.  Now, a singular expansion
compatible with~\eqref{eq:40} is given by
\begin{equation*}
  r^{[1]}(z) \sim \frac{-a_1}{2\sqrt\pi} \ln Z^{-1} + C_{1/2} + (Z\log
  Z)\A + Z\A, 
\end{equation*}
where the constant term is
\begin{equation}
  \label{eq:53}
  C_{1/2} := \sum_{n=m-1}^{\infty} \left( n^{1/2} \rho^n \tau_n
  + \frac{a_1}{2\sqrt\pi n} \right).
\end{equation}
Thus
\begin{equation*}
  r^{[1]}(z) + \sum_{j=0}^{m-2} x_j z^j \sim \frac{-a_1}{2\sqrt\pi} \ln
  Z^{-1} +   C_{1/2}' + (Z \log Z)\A + Z\A,
\end{equation*}
where
\begin{equation}
  \label{eq:54}
  C_{1/2}' := C_{1/2} + \sum_{j=0}^{m-2} x_j \rho^j.
\end{equation}
Using~\eqref{eq:3.7} and~\eqref{eq:5}, we get
\begin{equation}\label{eq:43}
  \begin{split}
  \muz{1} \sim \frac{a_0}{2\sqrt\pi(m-1)} Z^{-1/2}\ln Z^{-1} +
  \frac{a_0  C_{1/2}'}{-a_1(m-1)}Z^{-1/2} \\
  {} + (\log Z)\A + \A + (Z^{1/2} \log  Z)\A +   Z^{1/2}\A.
\end{split}
\end{equation}
Using singularity analysis and~\eqref{eq:13} we conclude
\begin{equation}
  \label{eq:42}
  \mun{1} \sim \frac{a_0}{-a_1\sqrt\pi(m-1)} n \ln{n} + \eta_{1/2} n
  + n^{1/2}\N + (\log n) \N + \N,
\end{equation}
where
\begin{equation*}
  \eta_{1/2} := \frac{2\sqrt\pi}{-a_1} \left( \frac{a_0(\gamma + 2
  \ln2)}{2\pi(m-1)} + \frac{a_0 C_{1/2}'}{-a_1\sqrt\pi(m-1)}
  \right).
\end{equation*}

\subsection{Higher moments}
\label{sec:higher-moments}

We will use induction to obtain asymptotics for higher-order moments.
Throughout $\alpha' := \alpha + \frac12$.  We consider the case
$\alpha > 1/2$ in Proposition~\ref{prop:alpha-gt-12} and handle the
remaining cases in Propositions~\ref{prop:alpha-lt-12}
and~\ref{prop:alpha=12}.

\begin{proposition}
  \label{prop:alpha-gt-12}
  Let $\alpha > 1/2$.  Then,  for $s \geq 1$, and $\epsilon > 0$ small enough,
  \begin{equation*}
    \muz{s} = D_s Z^{-s\alpha' + \frac12} + O(|Z|^{-s\alpha'+\frac12+q}),
  \end{equation*}
  where $q := \min\{\alpha - \frac12,\frac12\} - \epsilon$ with
  \begin{equation*}
    D_1 := \frac{a_0\Gamma(\alpha-\tfrac12)}{2(m-1)\sqrt\pi},
  \end{equation*}
  and, for $s \geq 2$,
  \begin{equation}
    \label{eq:8}
    D_s = \frac{a_0}{(m-1)(-a_1)} \left[ \frac{m-1}{2 a_0}
    \sum_{j=1}^{s-1} \binom{s}{j} D_j D_{s-j} +
    \frac{\Gamma(s\alpha'-1)}{\Gamma((s-1)\alpha'-\frac12)} s D_{s-1} \right].
  \end{equation}
\end{proposition}
\begin{proof}
  We proceed by induction on~$s$.  For $s=1$ the claim was proved
  as~\eqref{eq:10} and~\eqref{eq:39}.  [Note that $\muz{0} = \tau(z)
  \sim a_0$.]  Suppose $s \geq 2$.  We will first obtain the
  asymptotics of $r^{[s]}(z)$ at~\eqref{eq:2A.4} by analyzing each of
  the terms in the sum there.
  
  Suppose exactly $k \geq 1$ of $s_1,\dots,s_m$, say $s_1, \dots,
  s_k$, are nonzero.  Then, by induction,
  \begin{equation*}
    z^{m-1} \muz{s_1} \cdots \muz{s_m} = O(|Z|^{-(s-s_0)\alpha' +
    \frac{k}{2}}). 
  \end{equation*}
  Moreover, if $s_0=0$ then the contribution to $r^{[s]}(z)$ is
  $O(|Z|^{-s\alpha' + \frac32})$ unless $k=1$ or $k=2$.  (Observe,
  however, that if $k=1$ then $s_0$ cannot be zero as that would imply
  $s_1=s$.)  On the other hand, if $s_0 \ne 0$, then using singularity
  analysis for polylogarithms~\cite{MR2000a:05015} and Hadamard
  products~\cite{FFK}, we see that
  \begin{equation*}
    b^{\odot s_0}(z) \odot [ z^{m-1} \muz{s_1} \cdots \muz{s_m}  ] =
    O(|Z|^{-s\alpha' + \frac{s_0}2 + \frac{k}{2}}),
  \end{equation*}
  which is $O(|Z|^{-s\alpha' + \frac32 - \epsilon})$ unless $k=1$ and
  $s_0=1$.  (The $\epsilon$ term in the exponent avoids logarithmic
  factors that arise when $-s\alpha' + \frac{s_0}2 +
  \frac{k}{2}$ is a nonnegative integer.)

  If all of $s_1, \dots, s_m$ are zero, then $s_0=s$ and,
  using~\eqref{eq:11}, the contribution to $r^{[s]}(z)$ is
  $O(|Z|^{-s\alpha'+ \frac{s}2 + 
    \frac12})$ which is $O(|Z|^{-s\alpha' + \frac32})$.

  Hence unless $s_0=0$ and exactly two of $s_1,\dots,s_m$
  are nonzero or $s_0=1$ and exactly one of $s_1, \dots, s_m$ is
  $s-1$ in~\eqref{eq:2A.4}, the contribution to $r^{[s]}(z)$ is
  $O(|Z|^{-s\alpha' + \frac32 - \epsilon})$.  In the former case the
  contribution to $r^{[s]}(z)$ is gotten by using the induction
  hypothesis as
  \begin{equation*}
    \binom{m}2 \rho^{m-1} Z^{-s\alpha'+1} a_0^{m-2} \sum_{j=1}^{s-1}
    \binom{s}{j} D_{j} D_{s-j} + O(|Z|^{-s\alpha'+1+q}).
  \end{equation*}
  In the latter case, again using the induction hypothesis and
  singularity analysis for Hadamard products we get the contribution
  to $r^{[s]}(z)$ as
  \begin{equation*}
    m \rho^{m-1} a_0^{m-1} s D_{s-1}
    \frac{\Gamma(s\alpha'-1)}{\Gamma((s-1)\alpha'-\frac12)}
    Z^{-s\alpha'+1} + O(|Z|^{-s\alpha'+1+q}).
  \end{equation*}

  Finally, noting that the contribution from $\sum_{j=0}^{m-2} x_j^s
  z^j$  to the numerator on the right side in~\eqref{eq:3.7} is
  negligible, we complete the induction by using~\eqref{eq:4}
  and~\eqref{eq:5}.
\end{proof}

For $\alpha < 1/2$, it will be convenient to consider instead the
``approximately centered'' random variable
\begin{equation}\label{eq:44}
  \widetilde{X}_n := X_n - \frac{2a_0 C_{\alpha}}{(m-1)a_1^2} (n+1) =
  X_n - 
  \frac{\rho m^{\frac{m}{m-1}} C_\alpha}{(m-1)\alpha^*}(n+1),
\end{equation}
where $C_\alpha$ is defined at~\eqref{eq:41} and $\alpha^*$
at~\eqref{eq:18}.  See~\eqref{eq:9} for the motivation behind this
definition.  The choice of centering by a multiple of $n+1$ rather
than $n$ is motivated by the fact that with this centering
$\widetilde{X}_n$ satisfies the same distributional
recurrence~\eqref{eq:1.1} as $X_n$ [with appropriate initial
conditions $(\widetilde{X}_0, \ldots, \widetilde{X}_{m-2})$].  We will
use $\tilde{r}^{[s]}(z)$, $\tilde{\mu}^{[s]}(z)$, and
$\tilde{\mu}_n^{[s]}$ to denote the analogous quantities
for~$(\widetilde{X}_n)$.
\begin{proposition}
  \label{prop:alpha-lt-12}
    Let $\alpha < 1/2$.  Then,  for $s \geq 1$, and $\epsilon > 0$
    small enough, 
    \begin{equation*}
      \tilde{\mu}^{[s]}(z) = D_s Z^{-s\alpha' + \frac12} +
      O(|Z|^{-s\alpha'+1-\epsilon}) + c_s,
    \end{equation*}
    where $c_s$ is a constant and
    $D_s$ is defined as in Proposition~\ref{prop:alpha-gt-12}.
\end{proposition}
\begin{proof}[Proof outline]
  The basis of the induction is~\eqref{eq:52} and the induction step is
  identical to the one in the proof of
  Proposition~\ref{prop:alpha-gt-12}.  We omit the details.
\end{proof}

When $\alpha = 1/2$, we define
\begin{equation}\label{eq:45}
  \widetilde{X}_n := X_n - \frac{2a_0 C_{1/2}'}{(m-1) a_1^2} (n+1) =
  X_n - 
  \frac{\rho m^{\frac{m}{m-1}} C_{1/2}'}{(m-1)\alpha^*}(n+1).
\end{equation}
The constant~$C_{1/2}'$ is defined at~\eqref{eq:54}
using~\eqref{eq:53}.  The key result here is the following.
\begin{proposition}
  \label{prop:alpha=12}
  Let $\alpha = 1/2$.  Define $\sigma_m := -a_1(m-1)/(\sqrt2
  a_0)$.  Then, for $s \geq 1$,
  \begin{equation*}
    \tilde{\mu}^{[s]}(z) = -a_1 \sigma_m^{-s} Z^{-s+\frac12} \sum_{r=0}^s
    C_{s,r} (\ln Z^{-1})^{s-r} + O(|Z|^{-s + 1 - \epsilon}),
  \end{equation*}
  where the constants~$C_{s,r}$ do not depend on~$m$.
\end{proposition}
\begin{proof}[Proof sketch]
  The form of the proof is the same as those of
  Propositions~\ref{prop:alpha-gt-12} and~\ref{prop:alpha-lt-12}.  The
  basis of the induction is~\eqref{eq:43}.  [Note that without the
  centering we have done we would be saddled with the $Z^{-1/2}$ term,
  whose coefficient depends on~$m$, in~\eqref{eq:43}.]  For the
  induction step, in estimating $\tilde{r}^{[s]}(z)$, unless
  \begin{enumerate}
  \item $s_0=0$, and exactly two of $s_1, \dots, s_m$ are nonzero; or
  \label{item:5}
  \item $s_0=1$, and exactly one of $s_1, \dots, s_m$ is nonzero,
  \label{item:6}
  \end{enumerate}
  the contribution is $O(|Z|^{-s+\frac32-\epsilon})$.

  In case~(\ref{item:5}), by induction, the contribution to
  $\tilde{r}^{[s]}(z)$ is
  \begin{equation*}
    a_1^2 \sigma_m^{-s} \binom{m}{2} \rho^{m-1} a_0^{m-2} Z^{-s+1}
    \sum_{r=0}^s F_{s,r} (\ln Z^{-1})^{s-r} +
    O(|Z|^{-s+\frac32-\epsilon}),
  \end{equation*}
  where the~$(F_{s,r})$ do not depend on~$m$.  Using~\eqref{eq:4} and
  the definition of~$\sigma_m$, we see that the constant  multiplying
  the lead sum is
  \begin{equation*}
    a_1^2 \sigma_m^{-s} \binom{m}{2} \rho^{m-1} a_0^{m-2} = -a_1
    \sigma_m^{-(s-1)}/\sqrt2.
  \end{equation*}
  In case~(\ref{item:6}), the contribution is
  \begin{equation*}
    -a_1 \sigma_m^{-(s-1)} m \rho^{m-1} a_0^{m-1} Z^{-s+1}
     \sum_{r=0}^{s-1} G_{s,r} (\ln Z^{-1})^{s-1-r} +
     O(|Z|^{-s+\frac32-\epsilon}), 
  \end{equation*}
  where the $(G_{s,r})$ do not depend on~$m$.  Again,
  using~\eqref{eq:4}, the constant multiplying the lead sum is
  \begin{equation*}
    -a_1 \sigma_m^{-(s-1)} m \rho^{m-1} a_0^{m-1} = -a_1
     \sigma_m^{-(s-1)}.
  \end{equation*}
  Summing all the contributions we get
  \begin{equation*}
    \tilde{r}^{[s]}(z) = -a_1 \sigma_m^{-(s-1)} Z^{-s+1}
    \sum_{r=0}^{s} H_{s,r}(\ln Z^{-1})^{s-r} + O(|Z|^{-s+\frac32-\epsilon}),
  \end{equation*}
  where the $(H_{s,r})$ do not depend on~$m$.
  Recalling~\eqref{eq:3.7} and~\eqref{eq:5}, and the definition
  of~$\sigma_m$ once again, completes the induction.
\end{proof}
\begin{remark}\label{rem:m-free}
  The significance of Proposition~\ref{prop:alpha=12} is that the case
  $m=2$ has already been considered in~\cite{fill03:_limit_catal},
  allowing us to determine the desired limiting
  distribution~(Theorem~\ref{thm:alpha=half}) without 
  computing the constants $(C_{s,r})$ in
  Proposition~\ref{prop:alpha=12}.
\end{remark}

\subsection{Limiting distributions}
\label{sec:limit-distr}
We can now use the method of moments to derive limiting distributions
for the additive functional.

\begin{theorem}
  \label{thm:dist-ne-12}
  Let $\alpha \ne 1/2$, and let $X_n$ denote the additive functional that
  satisfies the distributional recurrence~\eqref{eq:1.1} with $b_n \equiv
  n^\alpha$.  Define $\alpha' := \alpha + \frac12$.
\begin{enumerate}[(a)]
\item If $\alpha > 1/2$, then
    \begin{equation*}
    (m-1)(m\alpha^*)^{1/2} \frac{X_n}{n^{\alpha'}} \toinlaw
    Y_{\alpha};
  \end{equation*}
\item if $\alpha < 1/2$, then
  \begin{equation*}
    \frac{(m-1)(m\alpha^*)^{1/2}}{n^{\alpha'}}\left[X_n -
    \frac{\rho m^{\frac{m}{m-1}} C_\alpha}{(m-1)\alpha^*}(n+1) \right] \toinlaw
    Y_{\alpha},
  \end{equation*}
  with $C_\alpha$ defined at~\eqref{eq:41} and $\alpha^*$ at~\eqref{eq:18}.
\end{enumerate}
In either case we have convergence of all moments, where $Y_{\alpha}$
has the unique distribution whose moments are given by $\E{Y_\alpha^s}
= M_s \equiv M_s(\alpha)$.  Here
  \begin{equation*}
    M_1 = \frac{\Gamma(\alpha-\frac12)}{\sqrt2\Gamma(\alpha)},
  \end{equation*}
  and, for $s \geq 2$,
  \begin{equation*}
    M_s = \frac{1}{4\sqrt\pi} \sum_{j=1}^{s-1} \binom{s}{j}
    \frac{\Gamma(j\alpha'-\frac12)
    \Gamma((s-j)\alpha'-\frac12)}{\Gamma(s\alpha'-\frac12)}
    M_j M_{s-j} + \frac{s\Gamma(s\alpha'-1)}{\sqrt2
    \Gamma(s\alpha'-\frac12)} M_{s-1}.
  \end{equation*}
\end{theorem}
\begin{proof}
  If $\alpha > 1/2$, then by Proposition~\ref{prop:alpha-gt-12},
  singularity analysis, and the asymptotics of $\tau_n$
  at~\eqref{eq:13}, we have
  \begin{equation*}
    \E X_n^s = \mun{s} = \frac{D_s
    2\sqrt\pi}{(-a_1)\Gamma(s\alpha'-\frac12)} n^{s\alpha'} +
    O(n^{s\alpha'-q}).
  \end{equation*}
  Define $\sigma \equiv \sigma_m := -a_1(m-1)/(\sqrt2 a_0) =
  (m-1)(\alpha^*/m)^{1/2}$, where the last equality uses~\eqref{eq:4},
  \eqref{eq:50}, and~\eqref{eq:18}.  Then, for fixed~$m$, as $n \to \infty$,
  \begin{equation*}
    \E{\left[  \sigma_m \frac{X_n}{n^{\alpha'}}
      \right]^s} \to M_s,
  \end{equation*}
  where, for $s \geq 1$,
  \begin{equation*}
    M_s := \frac{\sigma^s D_s 2
    \sqrt\pi}{(-a_1)\Gamma(s\alpha'-\frac12)}.
  \end{equation*}
  In particular, $M_1 = \Gamma(\alpha-\tfrac12)/[\sqrt2\Gamma(\alpha)]$.
  Furthermore, using~\eqref{eq:8}, we obtain the recurrence for $M_s$.

  Convergence in distribution follows from the fact that $(M_s)$
  satisfies Carleman's condition, as has been established
  in~\cite{fill03:_limit_catal}.

  The same proof holds if $\alpha < 1/2$, now by considering
  $\widetilde{X}_n$ defined at~\eqref{eq:44} and using
  Proposition~\ref{prop:alpha-lt-12}.
\end{proof}

The case~$\alpha = 1/2$ is covered by the following result.  As
alluded to in Remark~\ref{rem:m-free}, we will use the known results
for $m=2$ to derive the distribution for all~$m$.
\begin{theorem}
  \label{thm:alpha=half}
  Let $X_n$ denote the additive functional that satisfies the
  distributional recurrence~\eqref{eq:1.1} with $b_n \equiv n^{1/2}$.
  Define
  \begin{equation*}
    d_0(m) := \frac{\rho  m^{\frac{m}{m-1}} C_{1/2}'}{(m-1)\alpha^*},
  \end{equation*}
  where $C_{1/2}'$ is defined at~\eqref{eq:54} using~\eqref{eq:53}.
  Then, as $n \to \infty$,
  \begin{equation*}
    n^{-1} \left\{ \sigma_m X_n - \frac1{\sqrt{2\pi}} \left(n \ln n +
    \left[(2 \ln 2 + \gamma) + \sqrt{2\pi}\sigma_m d_0(m)\right]n \right)
    \right\} \toinlaw Y_{1/2},
  \end{equation*}
  with convergence of all moments, where $Y_{1/2}$ has the unique
  distribution with moments $m_k := \E{Y_{1/2}^k}$ given by $m_0 = 1$, $m_1
  = 0$, and for $k \geq 2$,
  \begin{multline*}
    m_k = \frac{1}{4\sqrt\pi} \frac{\Gamma(k-1)}{\Gamma(k-\frac12)} \\
    \times \left[ \sum_{\substack{k_1 + k_2 + k_3 = k\\ k_1, k_2 < k}}
    \binom{k}{k_1, k_2, k_3} m_{k_1} m_{k_2} \left(
    \frac1{\sqrt{2\pi}} \right)^{k_3} J_{k_1, k_2, k_3} + 4
    \sqrt\frac{\pi}{2} k m_{k-1} \right].
  \end{multline*}
  Here
  \begin{equation*}
    J_{k_1, k_2, k_3} := \int_0^1 x^{k_1 - \frac32} (1-x)^{k_2 -
    \frac32} \left[ x \ln x + (1-x) \ln{(1-x)} \right]^{k_3} \, dx.
  \end{equation*}
\end{theorem}
\begin{proof}
  Recall the definition of $\widetilde{X}_n \equiv \widetilde{X}_n(m)$
  at~\eqref{eq:45}.  Using
  Proposition~\ref{prop:alpha=12}, singularity analysis,
  and~\eqref{eq:13}, we see that
  \begin{equation*}
    \tmun{s} = \sigma_m^{-s} n^s \sum_{r=0}^s \widehat{C}_{s,r} (\log
    n)^{s-r} + O(n^{s-\frac12+\epsilon}),
  \end{equation*}
  where the $(\widehat{C}_{s,r})$ do not depend on~$m$.  Thus, for all
  (integer)~$s \geq 0$,
  \begin{equation*}
    \E{ [\sigma_m \widetilde{X}_n(m)]^s } = \E{[\sigma_2
    \widetilde{X}_n(2)]^s} +    O(n^{s-\frac12+\epsilon}).
  \end{equation*}
  It follows that
  \begin{align*}
   & \E{\left\{ \sigma_m X_n - \frac1{\sqrt{2\pi}} \left(n \ln n +
         \left[(2 \ln 2 + \gamma) + \sqrt{2\pi}\sigma_m d_0(m)\right]n
         \right) \right\}}^s \\
   &= \E{\left\{
         \sigma_m\widetilde{X}_n(m) - \frac1{\sqrt{2\pi}} [n \ln n + (2
         \ln 2 + \gamma)n]
       \right\}}^s \\
     &= \E{ \left\{ \sigma_2\widetilde{X}_n(2) - \frac1{\sqrt{2\pi}}
         [n \ln n + (2 \ln 2 + \gamma)n] \right\}}^s +
     O(n^{s-\frac12+\epsilon}) \\
     &= 2^{-s/2}\E{\left\{X_n(2) - \frac1{\sqrt\pi} n\ln n - D_1 n\right\}^s} +
         O(n^{s-\frac12 + \epsilon}),
  \end{align*}
  the last equality using $\sigma_2 = 2^{-1/2}$ and
  \begin{equation*}
    D_1 := \frac1{\sqrt\pi}[ 2 \ln 2 + \gamma + \sqrt\pi d_0(2)].
  \end{equation*}
  Observe that
  \begin{equation*}
    -\delta_n := \left( \frac1{\sqrt\pi} n \ln n + D_1 n \right) -
     \left[ \frac1{\sqrt\pi} (n+1) \ln{(n+1)} + D_1 (n+1) \right] =
     O(\log n).
  \end{equation*}
  Hence, for $l \geq 0$,
  \begin{align*}
    & \E{\left[X_n(2) - \frac1{\sqrt\pi} n\ln n - D_1 n\right]^s} \\
    &= \E{\left[X_n(2) - \frac1{\sqrt\pi} (n+1)\ln{(n+1)} - D_1 (n+1)
        + \delta_n \right]^s} \\
    &= \sum_{k=0}^s \binom{s}{r} [ m_k n^k + o(n^k) ] O( (\log
    n)^{s-k} ) \\
    & = [m_s + o(1)] n^s,
  \end{align*}
  where
  \begin{equation*}
  m_k := \lim_{n \to \infty} \E{\left\{n^{-1}\left[X_n(2) -
        \frac1{\sqrt\pi} 
        (n+1)\ln{(n+1)} -  D_1 (n+1) \right] \right\}^k }.
\end{equation*}
But the $m_k$'s have already been
determined~\cite[Proposition~3.8]{fill03:_limit_catal} and the claim
follows from there.
\end{proof}

\section{The shape functional}
\label{sec:shape-functional}

Recall that the shape functional~$X_n$ is the additive functional on $m$-ary
search trees induced by the toll
\[
b_n = \ln \binom{n}{m-1},
\]
with $(X_0,\ldots,X_{m-2}) = \mathbf{0}$. 
By~\eqref{eq:3.7}, we have
\begin{equation}
  \label{eq:16}
  \muz{s} = \frac{r^{[s]}(z)}{1 - m[z \tau(z)]^{m-1}},
\end{equation}
with $r^{[s]}(z)$ given by~\eqref{eq:2A.4}, where
\begin{equation*}
  b(z) = \sum_{n=m-1}^\infty \ln \binom{n}{m-1} z^n.
\end{equation*}
It follows from Theorem~1 of~\cite{MR2000a:05015} that $b(z)$ is amenable
to singularity analysis.
In the sequel we will make use of the following asymptotic expansion
of $b_n$ as $n \to \infty$:
\begin{equation}
  \label{eq:12}
  b_n \sim (m-1) \ln n - \ln[(m-1)!] + n^{-1}\N.
\end{equation}

\subsection{Mean}
\label{sec:mean-1}
\newcommand{\clog}{C_{\text{ln}}}

Using~\eqref{eq:14} and~\eqref{eq:12} we have
\begin{equation*}
  \begin{split}
  [z^n] r^{[1]}(z) \sim n^{-3/2} \rho^{-n}
  \biggl[
    \frac{(-a_1)(m-1)}{2\sqrt\pi} \ln n - \ln[(m-1)!]
    \left(\frac{-a_1}{2\sqrt\pi}\right) \\
    {}+ (n^{-1} \ln n) \N + n^{-1} \N
    \biggr].
  \end{split}
\end{equation*}
A compatible singular expansion can be computed using~\eqref{eq:15}:
\begin{align*}
  r^{[1]}(z) &\sim \clog - (-a_1)(m-1) Z^{1/2} \ln Z^{-1} \\
  & \quad {}- \left\{(-a_1) (m-1)[2(1 - \ln 2) - \gamma] -
    (-a_1)\ln[(m-1)!]
  \right\} Z^{1/2} \\
  & \quad {}+ Z\A + (Z^{3/2} \log Z)\A + Z^{3/2}\A,
\end{align*}
where
\begin{equation*}
  \clog := r^{[1]}(\rho) = \sum_{n=m-1}^\infty \rho^n \left[\ln
  \binom{n}{m-1}\right] \tau_n.
\end{equation*}
Now using~\eqref{eq:16} and~\eqref{eq:5} we get
\begin{equation}
  \label{eq:60}
  \begin{split}
  \muz{1} &\sim \frac{a_0\clog}{(-a_1)(m-1)}
  Z^{-1/2} \\
  & \quad {} - a_0 \ln Z^{-1} - a_0
  \left[
    \left(
      2(1 - \ln 2) - \gamma
    \right) - \frac{\ln[(m-1)!]}{m-1}
  \right]  + \frac{(m-2)\clog}{3(m-1)} \\
  & \qquad {} + (Z^{1/2} \log Z)\A  +
  Z^{1/2}\A + Z\A + (Z \log Z)\A,
\end{split}
\end{equation}
whence singularity analysis and~\eqref{eq:13} yield
\begin{equation*}
  \mun{1} \sim \frac{2 a_0 \clog}{(m-1)a_1^2} n - 2 \sqrt\pi
  \left(\frac{a_0}{-a_1}\right) n^{1/2} + (\log n)\N + \N + n^{-1/2}\N.
\end{equation*}

\subsection{Second moment and variance}
\label{sec:variance}

As in the case of the toll $n^{\alpha}$ when $\alpha < 1/2$, it will
be convenient to consider the random variable
\begin{equation*}
  \widetilde{X}_n := X_n - d_1(n+1),
\end{equation*}
where here
\begin{equation*}
  d_1 := \frac{2 a_0 \clog}{(m-1)a_1^2}.
\end{equation*}
Thus, $\widetilde{X}_n$ satisfies the same distributional
recurrence~\eqref{eq:1.1} as $X_n$ with initial conditions
$(\widetilde{X}_0, \ldots, \widetilde{X}_{m-2}) = -d_1(1,\dots,m-1)$.
Again, we use $\tilde{r}^{[s]}(z)$, $\tilde{\mu}^{[s]}(z)$, and
$\tilde{\mu}_n^{[s]}$ to denote the analogous quantities
for~$\widetilde{X}_n$.  Then, noting that $\E{\widetilde{X}_n} =
\E{X_n} - d_1(n+1)$, a singular expansion for $\tilde{\mu}^{[1]}(z)$
can be obtained using~\eqref{eq:60}, namely,
\begin{equation}
  \label{eq:17}
  \tilde{\mu}^{[1]}(z) = -a_0 \ln Z^{-1} - d_2 + O(|Z|^{\frac12 -
  \epsilon}),
\end{equation}
where
\begin{equation*}
  d_2 := a_0
  \left[
    \left(
      2(1 - \ln 2) - \gamma
    \right) - \frac{\ln[(m-1)!]}{m-1}
  \right]  - \frac{(m-2)\clog}{3(m-1)} + (a_0-a_2)d_1
\end{equation*}

We begin the variance computation by obtaining asymptotics for the
second moment~$\tilde{\mu}_n^{[2]}$.  To that end, we calculate the
contribution of the terms in the sum
\begin{equation*}
  \tilde{r}^{[2]}(z) = \sum_{\substack{s_0 + \dots + s_m = 2\\ s_1,
  \dots, s_m < 2}} \binom{2}{s_0, \dots, s_m} b^{\odot s_0}(z) \odot
\left[
  z^{m-1} \tilde{\mu}^{[s_1]}(z) \cdots \tilde{\mu}^{[s_m]}(z)
\right].
\end{equation*}

When $s_0=0$, exactly two of $s_1, \dots, s_m$ equal~1.  The
contribution to~$\tilde{r}^{[2]}(z)$ from such terms is
\begin{multline*}
  2 \binom{m}{2} z^{m-1} \tau^{m-2}(z)
  \left[\tilde{\mu}^{[1]}(z)\right]^2 \\
  =  m(m-1)\rho^{m-1} a_0^{m-2} \left[
  a_0^2 \ln^2 Z^{-1} + 2 a_0 d_2 \ln Z^{-1}  + d_2^2\right] + 
  O(|Z|^{\frac12-2\epsilon}).
\end{multline*}
When $s_0=1$, exactly one of $s_1, \dots, s_m$ equals~1.  The
contribution to~$\tilde{r}^{[2]}(z)$ from such terms is obtained
(after some routine calculations using the Zigzag algorithm) using the
expansion for $\Li_{1,1}$ at~\eqref{eq:15} as
\begin{equation*}
  \begin{split}
  & 2 m b(z) \odot \left[ z^{m-1} \tau^{m-1}(z) \tilde{\mu}^{[1]}(z)
  \right] \\
  & = 2 m \rho^{m-1} a_0^{m}
  \left[
    -\frac{m-1}2 \ln^2 Z^{-1} + (\gamma(m-1) + \ln[(m-1)!]) \ln
  Z^{-1}
  \right] \\
  & \qquad {}+ \text{constant} + O(|Z|^{\frac12 -\epsilon}).
\end{split}
\end{equation*}

Finally, when $s_0=2$ the contribution to $r^{[1]}(z)$ is $b^{\odot
  2}(z) \odot \left[z^{m-1} \tau^m(z) \right]$, which equals
  $\text{constant} + O(|Z|^{\frac12-\epsilon})$.

Summing these contributions and using~\eqref{eq:17} and~\eqref{eq:4},
we conclude (note the cancellation of the ostensible lead term) that
\begin{equation*}
  \tilde{r}^{[2]}(z) = 4 a_0(m-1)(1 - \ln 2) \ln Z^{-1} + d_3 +
  O(|Z|^{\frac12 - \epsilon}),
\end{equation*}
where 
\begin{equation*}
  d_3 := \lim_{z \to \rho} \left[r^{[2]}(z) - 4 a_0(m-1)(1 - \ln 2)
  \ln Z^{-1}\right].
\end{equation*}
This leads, using~\eqref{eq:5}, to
\begin{equation}
  \label{eq:19}
  \tilde{\mu}^{[2]}(z) = \frac{4 a_0^2}{-a_1} (1 - \ln 2)Z^{-1/2}
  \ln Z^{-1}
  + \frac{d_3 a_0}{(-a_1)(m-1)} Z^{-1/2} + O(|Z|^{-\epsilon}).
\end{equation}
By singularity analysis
\begin{equation*}
  \rho^n \tilde{\mu}_n^{[2]} \tau_n = \frac{4 a_0^2}{-\sqrt\pi a_1} (1
  - \ln 2)
  n^{-1/2} \ln n + d_4 n^{-1/2} + O(n^{-1+\epsilon}),
\end{equation*}
where
\begin{equation*}
  d_4 := \frac{4 a_0^2}{-\sqrt\pi a_1}(1-\ln2)(\gamma + 2 \ln 2) + \frac{d_3
  a_0}{\sqrt\pi(-a_1)(m-1)}.
\end{equation*}
Using the asymptotics of $\tau_n$ at~\eqref{eq:13} we get
\begin{equation*}
  \tilde{\mu}_n^{[2]} = 8 \left(\frac{a_0}{a_1}\right)^2 (1 - \ln 2) n
  \ln n + \frac{2 \sqrt \pi d_4}{-a_1} n + O(n^{\frac12+\epsilon}).
\end{equation*}
Thus
\begin{multline*}
  \Var{X_n} = \Var{\widetilde{X}_n} = \tilde{\mu}_n^{[2]} -
  \left(\tilde{\mu}_n^{[1]}\right)^2 \\
  = 8
  \left(\frac{a_0}{a_1}\right)^2 (1 - \ln 2) n \ln n + \left( \frac{2
  \sqrt  \pi d_4}{-a_1} - \frac{4 \pi a_0^2}{a_1^2} \right) n +
  O(n^{\frac12 + \epsilon}).
\end{multline*}

\subsection{Higher moments and limiting distribution}
\label{sec:limit-distr-1}

\begin{proposition}
  \label{prop:shape-moments}
  For $s \geq 2$ and $\epsilon > 0$ small enough,
  \begin{equation*}
    \tilde{\mu}^{[s]}(z) = Z^{-\frac{s}{2}+\frac12}
    \sum_{j=0}^{\floor{{s}/2}} C_{s,j} \left(\ln^{\floor{{s}/{2}} -
    j} Z^{-1} \right) + O(|Z|^{-\frac{s}{2} + 1 - \epsilon}),
  \end{equation*}
  with
  \begin{equation}
    \label{eq:24}
    C_{2l,0} = \frac{1}{-2a_1} \sum_{j=1}^{l-1} \binom{2l}{2j}
    C_{2j,0} C_{2l-2j,0}, \quad l \geq 2; \qquad C_{2,0} = \frac{4
    a_0^2}{-a_1} (1 - \ln 2).
  \end{equation}
\end{proposition}
\begin{proof}[Proof sketch]
  We proceed by induction.  For $s=2$ the claim is true
  by~\eqref{eq:19}.  [Recalling~\eqref{eq:17}, we note in passing that
  the claim is \emph{not} true for $s=1$.]  In the rest of the proof
  we will explicitly compute $C_{s,j}$ only when $s$ is even and
  $j=0$.  The rest of the coefficients will appear as
  unspecified constants.

For the induction step consider $s \geq 3$.  In a manner analogous to
the proof of Proposition~\ref{prop:alpha-gt-12} we obtain
the asymptotics of~$\tilde{r}^{[s]}(z)$ by first analyzing the
contributions of the terms in the sum~\eqref{eq:2A.4}.  As in that
proof one can check that unless
\begin{enumerate}[(a)]
\item $s_0=0$ and exactly two of $s_1, \dots, s_m$ are nonzero, or
  \label{item:1}
\item $s_0=1$ and exactly one of $s_1, \dots, s_m$ is $s-1$ \label{item:2}
\end{enumerate}%
in~\eqref{eq:2A.4}, the contribution to $\tilde{r}^{[s]}(z)$ is
\mbox{$O(|Z|^{-\frac{s}2 + \frac32 - (s+1)\epsilon})$}.

In case~(\ref{item:1}), the contribution to $\tilde{r}^{[s]}(z)$ is
\begin{equation}
  \label{eq:20}
  \binom{m}{2} z^{m-1} \tau^{m-2}(z) \sum_{j=1}^{s-1} \binom{s}{j}
  \tilde{\mu}^{[j]}(z) \tilde{\mu}^{[s-j]}(z).
\end{equation}
In this sum unless $j$ is 1 or $s-1$, the induction
hypothesis implies a contribution of
\begin{equation}
  \label{eq:21}
  \binom{m}2 \rho^{m-1}a_0^{m-2} Z^{-\frac{s}2 + 1}
  \sum_{l=0}^{\floor{{j}/2} + \floor{\frac{s-j}{2}}} \binom{s}{j}
  A_{s,j,l} \left(
    \ln^{\floor{{j}/2} + \floor{\frac{s-j}{2}} - l} Z^{-1}\right)
  +  O(|Z|^{-\frac{s}2 + \frac32 - 2\epsilon}),
\end{equation}
with $A_{s,j,0} = C_{j,0} C_{s-j,0}$.  Notice that when $s$ is even
and $j$ is odd, this contribution is $O\left(|Z|^{-\frac{s}{2}+1}
  \left(\ln^{\floor{{s}/2}-1} Z^{-1}\right)\right)$.  In all other
parity cases the contribution is {$O\left(|Z|^{-\frac{s}{2}+1}
    \ln^{\floor{{s}/2}} Z^{-1}\right)$}.

On the other hand the total contribution to $\tilde{r}^{[s]}(z)$ from
the terms where $j$ is 1 or $s-1$ in the sum in~\eqref{eq:20} is
obtained using the induction hypothesis and~\eqref{eq:17} as
\begin{equation}
  \label{eq:22}
  - 2s\binom{m}2 \rho^{m-1}a_0^{m-1} Z^{-\frac{s}2 + 1}
    \sum_{j=0}^{\floor{\frac{s-1}2}} E_{s-1,j}
    \ln^{\floor{\frac{s-1}2}+1-j} Z^{-1} + O(|Z|^{-\frac{s}{2} +
    \frac32 - 2\epsilon}),
\end{equation}
where $E_{s-1,0} = C_{s-1,0}$.

In case~(\ref{item:2}), the total contribution to $\tilde{r}^{[s]}(z)$
is
\begin{equation*}
  m s b(z) \odot [z^{m-1} \tau^{m-1}(z) \tilde{\mu}^{[s-1]}(z)].
\end{equation*}
Now another application of the Zigzag algorithm yields this
contribution to be
\begin{equation}
  \label{eq:23}
  s m(m-1)\rho^{m-1} a_0^{m-1} Z^{-\frac{s}2 + 1}
  \sum_{j=0}^{\floor{\frac{s-1}{2}}} D_{s,j}
  \ln^{\floor{\frac{s-1}2} + 1 -j} Z^{-1} + O(|Z|^{-\frac{s}{2} +
  \frac32 - 2\epsilon}),
\end{equation}
where $D_{s,0} = C_{s-1,0}$.  Note that the lead term here is exactly
the same as in~\eqref{eq:22} but with opposite sign, so
that~\eqref{eq:22} and~\eqref{eq:23} cancel each other to lead order.

Summing the various contributions, we find that $\tilde{r}^{[s]}(z)$ is
of the form
\[
\tilde{r}^{[s]}(z) = Z^{-\frac{s}2+1} \sum_{j=0}^{\lfloor s/2 \rfloor}
\widehat{C}_{s,j} \left( \ln^{\lfloor s/2 \rfloor -j} Z^{-1} \right) +
O\left( |Z|^{-\frac{s}{2} + \frac32 - 2\epsilon}\right),
\]
where, for $s$ even,
\[
\widehat{C}_{s,0} = \binom{m}2 \rho^{m-1} a_0^{m-2} \sum_{\substack{0
    < j < s\\\text{$j$ even}}} \binom{s}j C_{j,0} C_{s-j,0}.
\]
Using~\eqref{eq:5}, the result follows.
\end{proof}
It follows from Proposition~\ref{prop:shape-moments}, singularity
analysis, and the asymptotics of~$\tau_n$ at~\eqref{eq:13} that
\begin{equation*}
  \tilde{\mu}_n^{[s]} = \frac{2 \sqrt\pi
  C_{s,0}}{-a_1\Gamma(\frac{s-1}2)}n^{{s}/{2}}
  \ln^{\floor{{s}/{2}}}n + O(n^{{s}/2} \ln^{\floor{s/2}-1}n),
\end{equation*}
whence for $s \geq 1$, as $n \to \infty$,
\begin{equation*}
  \E{\left[ \frac{\widetilde{X}_n}{\sqrt{n\ln n}} \right]^{2s}} \to
  \frac{2\sqrt\pi C_{2s,0}}{(-a_1)\Gamma(s-\frac12)} \quad \text{and} \quad
  \E{\left[ \frac{\widetilde{X}_n}{\sqrt{n\ln n}} \right]^{2s-1}} = o(1).
\end{equation*}
Solving the recurrence~\eqref{eq:24} for $C_{2s,0}$ yields
\begin{equation*}
  C_{2s,0} = \frac{-a_1 \Gamma(\frac{s-1}2)}{2 \sqrt\pi}
  \frac{(2s)!}{2^s s!} \sigma^{2s} = \frac{-a_1}2 \frac{(2s)!
  (2s-2)!}{2^s 2^{2s-2}s!(s-1)!} \sigma^{2s},
\end{equation*}
where $\sigma^2 := 8 \left(a_0/a_1\right)^2 (1-\ln2)$.
The method of moments (see, for
example,~\cite[Theorem~30.1]{MR95k:60001}) implies then that the
shape functional is asymptotically normal.
\begin{theorem}
  \label{thm:shape-clt}
  Let $X_n$ denote the shape functional for uniformly distributed
  $m$-ary search trees on $n$~keys.  Then
  \begin{equation*}
    \frac{X_n - d_1(n+1)}{\sqrt{n \ln n}} \toinlaw
    {N}(0,\sigma^2) \quad \text{and} \quad \frac{X_n -
    \E{X_n}}{\sqrt{\Var{X_n}}} \toinlaw {N}(0,1),
  \end{equation*}
  where
  \begin{equation*}
    d_1 := \frac{2 a_0}{(m-1)a_1^2} \sum_{n=m-1}^\infty \rho^n \left[\ln
  \binom{n}{m-1}\right] \tau_n
  \end{equation*}
  and $\sigma^2 := 8 \left(a_0/a_1\right)^2 (1-\ln2)$.
\end{theorem}
\begin{remark}
  \label{rem:rp-vs-uniform}
  It is known~\cite{MR1871558,FK-transfer} that under the random
  permutation model the shape functional centered by its mean and
  scaled by its standard deviation is asymptotically normal for $2
  \leq m \leq 26$ and does not have a limiting distribution for $m >
  26$.  In contrast, under the uniform model we have asymptotic
  normality for all $m \geq 2$.
\end{remark}

\section{The space requirement}
\label{sec:space-requirement}

The space requirement for $m$-ary search trees is the number of
{nodes} in the tree~\cite{MR93f:68045}.  (For a binary search tree
the space requirement for $n$~keys is clearly~$n$, so in this section
we assume $m \geq 3$.)  The limiting
distribution of this parameter under the random permutation model has
been considered by several
authors~\cite{MR90a:68012,MR95i:68030,MR1871558,FK-transfer}.  In our
framework it is the additive functional~$X_n$ corresponding to the
toll~$\mathbf{1}_{n \geq m-1}$ with initial conditions \mbox{$(X_0,
  \ldots, X_{m-2}) = (0,1,\ldots,1)$}.

The $s$th moment $\muz{s} := \E{X_n^s}$ can be computed as usual
using~\eqref{eq:3.7}, 
where now
\begin{equation}
  \label{eq:25}
  r^{[s]}(z) = z^{m-1} \sum_{\substack{s_0 + \dots + s_m = s\\s_1,
  \dots, s_m < s}} \binom{s}{s_0, \dots, s_m} \muz{s_1} \cdots \muz{s_m}
\end{equation}
since the toll generating function \mbox{$b(z) = (1-z)^{-1}$} serves as
the identity for Hadamard products.

\subsection{Mean}
\label{sec:mean-2}
Substituting $s=1$ in~\eqref{eq:25} and using~\eqref{eq:3}
and~\eqref{eq:1} yields
\begin{equation*}
  r^{[1]}(z) + \sum_{j=1}^{m-2} z^j = z^{m-1}\tau^{m}(z) +
  \sum_{j=1}^{m-2} z^j = \tau(z) - 1 \sim (a_0-1) + a_1 Z^{1/2} + Z\A
  + Z^{3/2} \A.
\end{equation*}
Then, by~\eqref{eq:3.7} and~\eqref{eq:5},
\begin{equation}
  \label{eq:28}
    \muz{1} \sim \frac{a_0(a_0-1)}{-a_1(m-1)} Z^{-1/2} +
  \left[
    c_0(a_0-1) - \frac{a_0}{m-1}
  \right] + Z^{1/2}\A + Z\A.
\end{equation}
Singularity analysis and the asymptotics of~$\tau_n$ immediately lead
to
\begin{equation}\label{eq:35}
  \mun{1} \sim \frac{2 a_0(a_0-1)}{a_1^2 (m-1)} n + \N.
\end{equation}

\subsection{Variance}
\label{sec:variance-1}

As for the shape functional it is convenient to consider instead the
``centered'' functional $\widetilde{X}_n := X_n - d_1(n+1)$, where now
\begin{equation*}
  d_1 := \frac{2 a_0(a_0-1)}{a_1^2 (m-1)} = \frac{m(1-\rho
  m^\frac{1}{m-1})}{(m-1)\alpha^*}.
\end{equation*}
The ``centered'' space requirement satisfies the same distributional
recurrence as the space requirement with initial conditions
\[
(\widetilde{X}_0, \dots, \widetilde{X}_{m-2}) = -(d_1,
2d_1-1, \dots, (m-1)d_1-1).
\]
We employ the same notation as Section~\ref{sec:shape-functional},
with $\trz{s}$, $\tmuz{s}$, and~$\tmun{s}$ denoting quantities analogous
to $r^{[s]}(z)$, $\muz{s}$, and $\mun{s}$.

By definition
\begin{equation}\label{eq:29}
  \tmuz{1} = \muz{1} - d_1 \sum_{n=0}^\infty (n+1) \tau_n z^n.
\end{equation}
Now
\begin{equation}\label{eq:30}
  \sum_{n=0}^\infty (n+1) \tau_n z^n = z\tau'(z) + \tau(z) \sim
  \frac{-a_1}2 Z^{-1/2} + (a_0-a_2) + Z^{1/2}\A + Z\A.
\end{equation}
Thus using~\eqref{eq:28},~\eqref{eq:29}, and~\eqref{eq:30} we find
\begin{equation}\label{eq:31}
  \tmuz{1} \sim B_1 + Z^{1/2}\A + Z\A,
\end{equation}
where, using~\eqref{eq:49}, we have
\begin{equation}\label{eq:26}
  B_1 := c_0(a_0 - 1) - \frac{a_0}{m-1} - d_1 (a_0-a_2) =
  -\frac{a_0}{m-1}.
\end{equation}

Using~\eqref{eq:25},
\begin{equation}\label{eq:32}
  \trz{2} = z^{m-1}
  \left[
    m(m-1)
    \left(
      \tmuz{1}
    \right)^2 \tau^{m-2}(z) + 2m\tmuz{1} \tau^{m-1}(z) + \tau^{m}(z)
  \right].
\end{equation}
Also,
\begin{equation}\label{eq:33}
  \sum_{j=0}^{m-2} \tilde{x}_j^2 \rho^j = d_1^2 + \sum_{j=1}^{m-2}
  [-(j+1)d_1 + 1]^2 \rho^j =: \delta_1 .
\end{equation}
By~\eqref{eq:31}
\begin{equation*}
  \left(
    \tmuz{1}
  \right)^2 \sim B_1^2 + Z^{1/2}\A + Z\A,
\end{equation*}
and, for $k \geq 1$,
\begin{equation*}
  \tau^{k}(z) \sim a_0^k + k a_0^{k-1} a_1 Z^{1/2} + Z\A + Z^{3/2}\A,
\end{equation*}
so that
\begin{equation*}
  \left(
    \tmuz{1}
  \right)^2 \tau^{m-2}(z) \sim a_0^{m-2} B_1^2 + Z^{1/2}\A + Z\A.
\end{equation*}
Similarly
\begin{equation*}
  \tmuz{1} \tau^{m-1}(z) \sim a_0^{m-1} B_1 + Z^{1/2}\A + Z\A.
\end{equation*}
Using these expansions and~\eqref{eq:4} in~\eqref{eq:32} gives
[recalling~\eqref{eq:6} and~\eqref{eq:33}]
\begin{equation*}
  \trz{2} + \sum_{j=0}^{m-2} \tilde{x}_j^2 z^j \sim
  \left[
    \frac{m-1}{a_0} B_1^2 + 2 B_1 + \frac{a_0}{m} + \delta_1
  \right] + Z^{1/2}\A + Z\A,
\end{equation*}
whence
\begin{equation}
  \label{eq:34}
  \tmuz{2} \sim B_2 Z^{-1/2} + \A + Z^{1/2}\A,
\end{equation}
where
\begin{equation}\label{eq:27}
  B_2 := \frac{a_0}{-a_1(m-1)}
  \left[
    \frac{m-1}{a_0} B_1^2 + 2 B_1 + \frac{a_0}{m} + \delta_1
  \right] = \frac{a_0}{-a_1(m-1)} \left[ \delta_1 - \frac{a_0}{m(m-1)}
  \right].
\end{equation}
By singularity analysis and the asymptotics of~$\tau_n$, then
\begin{equation*}
  \tmun{2} \sim \frac{2 B_2}{-a_1} n + \N.
\end{equation*}
Recalling~\eqref{eq:35}, we observe that $\tmun{1} \sim \N$ so that
\begin{equation}
  \label{eq:36}
  \Var{X_n} = \Var{\widetilde{X}_n} = \tmun{2} - \left( \tmun{1}
  \right)^2 \sim
  \frac{2 B_2}{-a_1} n + \N.
\end{equation}
We pause here to remark that since $\Var{X_n} \to \infty$ as $n \to \infty$
for $m \geq 3$ (see Remark~\ref{rem:space-not-degenerate} in the
Appendix), we must have $B_2 > 0$.  We do not know a direct proof of
this fact.

\subsection{Higher moments and limiting distribution}
\label{sec:high-moments-limit}

\begin{proposition}
  \label{prop:space-moments}
  For $s \geq 1$,
  \begin{equation*}
    \tmuz{s} = B_s Z^{-\frac{s}2 + \frac12} + O(|Z|^{-\frac{s}{2} + 1}),
  \end{equation*}
  where $B_1$ and $B_2$ are given at~\eqref{eq:26} and~\eqref{eq:27},
  respectively, and, for $s \geq 3$,
  \begin{equation}
    \label{eq:38}
    B_s = \frac{a_0}{-a_1(m-1)}
    \left[
    \frac{m-1}{2 a_0} \sum_{j=1}^{s-1} \binom{s}{j} B_j B_{s-j} + s B_{s-1}
    \right].
  \end{equation}
\end{proposition}
\begin{proof}[Proof Sketch]
  We proceed by induction.  For $s=1$ and $s=2$, the result has been
  established at~\eqref{eq:31} and~\eqref{eq:34}, respectively.

  Suppose $s \geq 3$.  We will first obtain a singular expansion for
  $\trz{s}$ by analyzing the contributions of the terms in the sum
  at~\eqref{eq:25}.  As in the proofs of
  Propositions~\ref{prop:alpha-gt-12} and~\ref{prop:shape-moments}
  unless
  \begin{enumerate}[(a)]
  \item $s_0=0$ and exactly two of $s_1,\dots,s_m$ are nonzero,
  or \label{item:3}
  \item $s_0=1$ and exactly one of $s_1,\dots,s_m$ is $s-1$,
  \label{item:4}
  \end{enumerate}
  the contribution to $\tilde{r}^{[s]}(z)$ is
  $O(|Z|^{-\frac{s}2+\frac32})$. 

  In case~(\ref{item:3}), the contribution to $\trz{s}$ is
  \begin{equation*}
    \binom{m}2 \rho^{m-1} a_0^{m-2} Z^{-\frac{s}{2}
    + 1} \sum_{j=1}^{s-1} \binom{s}{j} B_j B_{s-j}  +
    O(|Z|^{-\frac{s}{2} + \frac32}).
  \end{equation*}
  In case~(\ref{item:4}), the contribution to $\trz{s}$ is
  \begin{equation*}
    m \rho^{m-1} s a_0^{m-1} B_{s-1} Z^{-\frac{s}{2}
    + 1} + O(|Z|^{-\frac{s}{2} + \frac32})
  \end{equation*}

This leads to
  \begin{equation*}
    \trz{s} + \sum_{j=0}^{m-2} \tilde{x}_j^s z^j =
    \left(
      \frac{m-1}{2 a_0} \sum_{j=1}^{s-1} \binom{s}{j} B_j B_{s-j} + s B_{s-1}
    \right) Z^{-\frac{s}2+1} + O(|Z|^{-\frac{s}{2} + \frac32}).
  \end{equation*}
whence~\eqref{eq:3.7} and~\eqref{eq:5} complete the induction.
\end{proof}
Now using~\eqref{eq:38} and~\eqref{eq:26} we have $B_3 = 0$ and by
induction $B_{2s+1} = 0$ for $s=1, 2, \dots$.  Then [compare~\eqref{eq:24}]
\begin{equation*}
  B_{2l} = \frac{1}{-2a_1} \sum_{j=1}^{l-1} \binom{2l}{2j} B_{2j}
  B_{2l-2j}, \quad l \geq 2,
\end{equation*}
with $B_2$ given at~\eqref{eq:27}.  Now following the development
leading to Theorem~\ref{thm:shape-clt} we can conclude asymptotic
normality for the space requirement.
\begin{theorem}
  \label{thm:space-clt}
  Let $X_n$ denote the space requirement for uniformly distributed
  $m$-ary search trees on $n$~keys.  Then
  \begin{equation*}
    \frac{X_n - d_1(n+1)}{\sqrt{n}} \toinlaw
    {N}(0,\sigma^2) \quad \text{and} \quad \frac{X_n -
    \E{X_n}}{\sqrt{\Var{X_n}}} \toinlaw {N}(0,1),
  \end{equation*}
  where
  \begin{equation*}
    d_1 = \frac{m(1-\rho m^\frac{1}{m-1})}{(m-1)\alpha^*}  \quad
    \text{and} \quad \sigma^2 = 2 \frac{B_2}{-a_1}.
  \end{equation*}
  Here $B_2$ is given by~\eqref{eq:27}, with $\delta_1$ defined
  at~\eqref{eq:33} and  $a_0$ and $a_1$ at~\eqref{eq:4}
  and~\eqref{eq:50}, respectively.
\end{theorem}

\section{Number of leaves}
\label{sec:number-leaves}

Lastly we consider the number of leaves in an $m$-ary search tree.
This is the additive functional~$X_n$ corresponding to the
toll~$\mathbf{1}_{n=m-1}$ with initial conditions~$(X_0, \dots,
X_{m-2}) = (0, 1, \dots, 1)$.  Under the random permutation model, the
number of leaves is asymptotically
normal~\cite[Theorem~2.5]{FK-transfer}.  We will establish an
analogous result (Theorem~\ref{thm:leaves-clt}) under the uniform model.

The toll generating function is now $b(z) = z^{m-1}$.  Note that
\begin{equation*}
  b^{\odot s}(z) = 
  \begin{cases}
    (1-z)^{-1} & s = 0, \\
    z^{m-1} & s \geq 1.
  \end{cases}
\end{equation*}
Thus in~\eqref{eq:2A.4}
\begin{equation}
  \label{eq:58}
  r^{[s]}(z) = z^{m-1} \left[ \sum_{\substack{s_1 + \dots + s_m = s\\
  s_1, \dots, s_m < s}} \binom{s}{s_1, \dots, s_m} \muz{s_1} \cdots
  \muz{s_m} +  1 \right].
\end{equation}
Given the similarity of the calculations with those of
Section~\ref{sec:space-requirement}, we will be brief.  The interested
reader is invited to flesh out the development of this section along
the lines of Section~\ref{sec:space-requirement}.

\subsection{Mean}
\label{sec:mean-3}

Substituting $s=1$ in~\eqref{eq:58} and using~\eqref{eq:3.7} we get
\begin{equation*}
  \muz{1} \sim \frac{a_0}{-a_1 m^{\frac{m}{m-1}}} Z^{-1/2} +
  \frac{m-2}{3 m^{\frac{m}{m-1}}} + Z^{1/2}\A + Z\A
\end{equation*}
so that
\begin{equation*}
  \mun{1} \sim \frac{2a_0}{a_1^2 m^{\frac{m}{m-1}}} n + \N \sim
  \frac{\rho}{\alpha^*} n + \N.
\end{equation*}

\subsection{Variance}
\label{sec:variance-2}

For the variance we consider again $\widetilde{X}_{n} := X_n -
d_1(n+1)$, where now $d_1 := \frac{\rho}{\alpha^*}$.  Then
\begin{equation*}
  \tmuz{1} \sim B_1 + Z^{1/2}\A + Z\A,
\end{equation*}
where [compare~\eqref{eq:26}] $B_1 = 0$.  Substituting $s=2$
in~\eqref{eq:58} we get
\begin{equation*}
  \tilde{r}^{[2]}(z) = z^{m-1}\left[ 1 + m(m-1) \left( \tmuz{1}
  \right)^2 \tau^{m-2}(z) \right] \sim \rho^{m-1} + Z\A + Z^{3/2}\A.
\end{equation*}
Thus
\begin{equation*}
  \tmuz{2} \sim B_2 Z^{-1/2} + \A + Z^{1/2}\A,
\end{equation*}
where
\begin{equation}
  \label{eq:59}
  B_2 = \frac{a_0(\rho^{m-1} + \delta_1)}{-a_1(m-1)},
\end{equation}
with $\delta_1$ defined as at~\eqref{eq:33} (but now with $d_1 =
\rho/\alpha^*$).  This leads to
\begin{equation*}
  \Var{X_n} \sim \frac{2B_2}{-a_1} n + \N \sim \frac{\rho
  m^{\frac{m}{m-1} }(\rho^{m-1} + \delta_1)}{\alpha^* (m-1)}n + \N.
\end{equation*}
Clearly $B_2 > 0$, and therefore $\Var{X_n}$ grows at an
asymptotically linear rate.

\subsection{Higher moments and limiting distribution}
\label{sec:high-moments-limit-1}

The analog of Proposition~\ref{prop:space-moments} is the following.
We omit the proof.
\begin{proposition}
  \label{prop:leaves-moments}
  For $s \geq 1$,
  \begin{equation*}
    \tmuz{a} = B_s Z^{-\frac{s}2 + \frac12} + O(|Z|^{-\frac{s}2 +
    1}),
  \end{equation*}
  where $B_1 = 0$, $B_2$ is given at \eqref{eq:59}, and, for $s \geq
  3$,
  \begin{equation*}
    B_s = - \frac{1}{-2a_1} \sum_{j=1}^{s-1} \binom{s}{j} B_j B_{s-j}.
  \end{equation*}
\end{proposition}
It follows easily from Proposition~\ref{prop:leaves-moments} that the
number of leaves is asymptotically normal.
\begin{theorem}
  \label{thm:leaves-clt}
  Let $X_n$ denote the number of leaves in a uniformly distributed
  $m$-ary search tree on $n$~keys.  Then
    \begin{equation*}
    \frac{X_n - \frac{\rho}{\alpha^*}(n+1)}{\sqrt{n}} \toinlaw
    {N}(0,\sigma^2) \quad \text{and} \quad \frac{X_n -
    \E{X_n}}{\sqrt{\Var{X_n}}} \toinlaw {N}(0,1),
  \end{equation*}
  where
  \begin{equation*}
    \sigma^2 = \frac{\rho
      m^{\frac{m}{m-1} }(\rho^{m-1} + \delta_1)}{\alpha^* (m-1)}. 
  \end{equation*}
\end{theorem}

\appendix

\section{Growth of variance for nondegenerate functionals}
\label{sec:degenerate}

In Section~\ref{sec:variance-1} we claimed that the variance of the
space requirement tends to infinity as $n \to \infty$ for $m \geq 3$.
In this appendix we identify additive functionals that are
degenerate, i.e., for each fixed~$n$ have the same value for all $m$-ary
search trees with $n$~keys, and provide a lower bound on the rate of
growth of the variance of any nondegenerate additive functional.
\begin{theorem}
  \label{thm:degenerate}
  Consider an additive functional $(X_n)_{n \geq 0}$ \emph{[}as
  at~\eqref{eq:1.1}\emph{]} with
  toll~$(b_n)_{n \geq 0}$, 
  with initial conditions $(x_0, \dots, x_{m-2}) = (b_0, \dots,
  b_{m-2})$.   Then the following statements are equivalent:
  \begin{enumerate}[(a)]
  \item $\mathcal{L}(X_n)$ is degenerate for every $n \geq
    0$. \label{item:7} 
  \item The toll satisfies
    \begin{equation*}
      b_n = 
      \begin{cases}
        n b_1 - (n-1) b_0, &  \quad n = 2, \dots, m-2 \\
        (m-1)(b_1 - 2b_0), &  \quad n \geq m-1.
      \end{cases}
    \end{equation*} \label{item:8}
  \item $X_n = n b_1 - (n-1) b_0$ for every $n \geq 0$. \label{item:9}
  \end{enumerate}
  Moreover, if~(\ref{item:7}) does not hold then $\sigma_n^2 :=
  \Var{X_n} = \Omega(\log n)$ as $n \to \infty$.
\end{theorem}
\begin{remark}
  \label{rem:space-not-degenerate}
  Before we prove the theorem we apply it to the space requirement.
  It is easily checked that, for this additive functional,
  condition~(\ref{item:8}) holds only when $m=2$. Thus the space
  requirement is not degenerate (and its variance tends to infinity
  as $n \to \infty$) for $m \geq 3$.
\end{remark}
\begin{proof}[Proof of Theorem~\ref{thm:degenerate}]
  We will show below that~(\ref{item:8}) and~(\ref{item:9}) are
  equivalent.

  To show the equivalence with~(\ref{item:7}), suppose first
  that~(\ref{item:7}) holds, so that $X_n = x_n$ deterministically for
  all $n \geq 0$.  Then
  \begin{align}
    \label{eq:51}
    x_{n+(m-1)} &= x_n + (m-1) x_0 + b_{n+(m-1)} && \text{for } n \geq
    0 \\
    &= x_{n-1} + x_1 + (m-2) x_0 + b_{n+(m-1)}   && \text{for } n \geq 1,
    \notag
  \end{align}
  and so $x_n = x_{n-1} + (x_1 - x_0) = x_{n-1} + (b_1-b_0)$ for $n
  \geq 1$.
  Condition~(\ref{item:9}) then follows by induction.  Conversely,
  condition~(\ref{item:9}) trivially implies~(\ref{item:7}), and the
  equivalent~(\ref{item:8}) shows that $X_n$ is indeed an additive
  functional.

  We now show the equivalence of~(\ref{item:8}) and~(\ref{item:9}).
  If~(\ref{item:9}), holds then so does~\eqref{eq:51} [because $(X_n) =
  (x_n)$ is an additive functional], from which [by solving for
  $b_{n+(m-1)}$] it is easy to check
  that~(\ref{item:8}) holds.  Conversely, if~(\ref{item:8}) holds,
  then~(\ref{item:9}) is trivially true for $n = 0, \dots, m-2$, and
  holds by induction for $n \geq m-1$.
  
  Suppose now that $(b_n)$ does not satisfy~(\ref{item:8}).  Let $n_0$
  be such that $\mathcal{L}(X_{n_0})$ is not degenerate. Then for any
  $n \geq n_1 := n_0 + (m-1)$, there is positive probability of having
  a subtree containing precisely $n_0$~keys, so that by the law of
  total variance we must $\sigma_n^2 > 0$.
  
  Finally we show that if $\mathcal{L}(X_n)$ is nondegenerate for
  all~$n \geq n_1$, then $\sigma_n^2 = \Omega(\log n)$ as $n \to
  \infty$.  First, it is clear that in this case there exists
  $\epsilon > 0$ such that $\sigma_n^2 \geq \epsilon$ for all $n \in
  [n_1, m n_1 + (m-2)]$.  Now suppose $n \geq m n_1 + (m-1)$.  Then at
  least one subtree must have size in the range $[n_1, n-(m-1)]
  \subseteq [n_1, n-1]$, so that by induction and the law of total
  variance we have $\sigma_n^2 \geq \epsilon$ for all $n \geq n_1$.

  For $n \geq n_1 + (m-1)$, let $p_n$ denote the probability that a
  tree of size~$n$ has its first subtree of size $n - n_1 - (m-1)$,
  its second of size~$n_1$, and the rest of size~0.  Then
  by~\eqref{eq:55} and the asymptotics of~$\tau_n$ (see
  Remark~\ref{rem:tau_n_asymptotics}), as $n \to \infty$ we have
  \begin{equation*}
    p_n = \frac{\tau_{n - n_1 - (m-1)} \tau_{n_1}}{\tau_n} \to
    \tau_{n_1} \rho^{n_1 + (m-1)} > 0.
  \end{equation*}
  Also $p_n > 0$ for $n \geq n_1 + (m-1)$, and so $\delta := \inf
  \{p_n: n
    \geq n_1 + (m-1) \} > 0$.
  
  Define $\alpha_0 := n_1$, $\alpha_1 := m n_1 + (m-1)$, and, for $k
  \geq 2$, $\alpha_k := m \alpha_{k-1} + n_1 + (m-1)$.  We will show,
  for each $k \geq 0$ and any $n \geq \alpha_k$,  that $\sigma_n^2
  \geq (1 + k \delta)\epsilon$, and the logarithmic-growth claim follows.
  
  For $k=0$ the result has been shown above.  For $k=1$, using the law
  of total variance and $p_n \geq \delta$ the result follows (compare
  the case $k \geq 2$ to follow). Suppose
  $k \geq 2$ and $n \geq \alpha_k$.  Then $n-(m-1) \geq m
  \alpha_{k-1}$, so there must be at least one subtree of size at
  least~$\alpha_{k-1}$.  Hence $\sigma_n^2 \geq (1 + (k-1)\delta)\epsilon$.
  But with probability~$p_{n} \geq \delta$, the first two subtrees are
  each of size at least~$n_1$; then by the law of total variance we have
  $\sigma_n^2 \geq [1 + (k-1)\delta]\epsilon + \delta \epsilon =
  (1+k\delta)\epsilon$, as desired.
\end{proof}

\bibliographystyle{habbrv}
\bibliography{master}

\bigskip

\end{document}